\documentclass[12pt]{article}
\usepackage{latexsym}
\usepackage{epsfig}
\usepackage{amsfonts,amssymb,amsmath}
\usepackage{graphicx,multicol,color}
\usepackage{stmaryrd}
\usepackage{float}
\usepackage{caption}
\title{An estimation of the stability and the localisability functions of multistable processes}
\author{R. Le Gu\'evel\\
\small{{\em Universit\'{e} de Rennes 2 - Haute Bretagne, Equipe de Statistique Irmar}}\\
\small{{\em Place du Recteur Henri Le Moal, CS 24307, 35043 RENNES cedex, France}}\\
\small{{\em ronan.leguevel@univ-rennes2.fr}}}

\date{}

\def\bbbr{{\bf R}} 
\def\bbbn{{\bf N}} 
\def\bbbz{{\bf Z}}
\newtheorem{theo}{Theorem}

\newtheorem{prop}[theo]{Proposition}

\newtheorem{lem}[theo]{Lemma}

\newcommand\E{\mbox{\sf E}}

\newcommand{\argmin}{\mathop{\mathrm{arg\,min}}}
\newcommand\cd{\stackrel{d}{\rightarrow}}

\renewcommand\P{{\sf P}}
\newcommand\Var{{\sf Var}}

\newcommand{\one}{\ifmmode {\sf 1}\hspace{-.26em}{\sf
l}\hspace{-.35em}{\sf \_} \else ${\sf 1}\hspace{-.26em}{\sf
l}\hspace{-.35em}{\sf \_}$ \fi}

\renewcommand{\Box}{\mbox{\rule{1ex}{1ex}}}

\begin{document}

\maketitle

\begin{abstract}
\noindent  Multistable processes are tangent at each point to a stable process, but where the index of stability and the index of localisability varies along the path. In this work, we give two estimators of the stability and the localisability functions, and we prove the consistency of those two estimators. We illustrate these convergences with two examples, the Levy multistable process and the Linear Multifractional Multistable Motion.
\end{abstract}

\vspace{1cm}


\vspace{1cm}

%

\section{Introduction}
\begin{figure}[H]
\begin{center}
\begin{tabular}{l}
      \includegraphics[height=4cm, width=15cm]{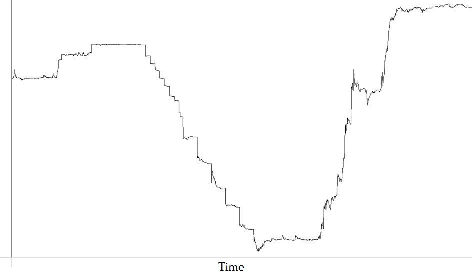} \\
\end{tabular}
\end{center}
\caption{Financial data where the increments do not appear to be stationary : the intensity of jumps is varying over time.}\label{trajFedFund}
\end{figure}

\begin{figure}[H]
\begin{center}
\begin{tabular}{ll}
      \includegraphics[height=4cm, width=13cm]{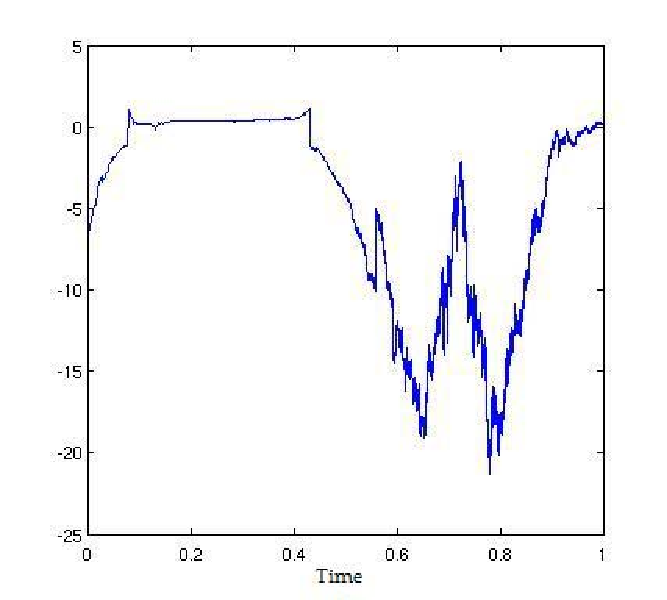} \\
\end{tabular}
\end{center}
\caption{Realization of a simulated multistable process.The sample size is $n=20000$.}\label{trajMulti}
\end{figure}

Multifractional multistable processes have been recently introduced as models for phenomena where the regularity and the intensity of jumps are non constant, and particularly when the increments of the observed trajectories are not stationary. 
In Figure \ref{trajFedFund}, we display a path of a financial data from federal funds, where the frequency of the jumps seems to vary with time.
The multistable processes then extend the stable models in order to take into account this additional variability (see Figure \ref{trajMulti} for an example of a realization of such a process, computed with the simulation method explained in \cite{FLGLV}).
We describe then some events with a low intensity of jumps at some times, which may be very erratic at other times. We provide another example of application in Figure \ref{EcgTraj} of Section \ref{AppliECG}, where we consider a path coming from electrocardiogram.

Multistable processes are stochastic processes which are locally stable, but where the index of stability $\alpha$ varies with ``time'', and therefore is a function.
They were constructed in \cite{FLGLV,FLV,FL,LGLV} using respectively moving averages, sums over Poisson processes, multistable measures, and the Ferguson-Klass-LePage series representation, this last definition being the representation used hereafter. 
These processes are, under general assumptions locally self-similar, with an index of self-similarity $H$ which is also a function.
In the remaining of this work, given one trajectory of a multistable process, we provide an estimator for each function. 

The aim of this work is then to introduce, for a large class of multistable processes, an estimator of the local index of stability $\alpha$. 
We prove in the sequel the consistency of this estimator with a convergence in all the $L^r$ spaces. This class includes two examples considered in \cite{FLV,LGLV}, the L\'evy multistable motion and linear multifractional 
multistable motion. We then estimate the local self-similarity function $H$. For the same class of multistable processes, we obtain a consistent estimator of $H$.
In the case of the L\'evy multistable motion, we are able to ascertain the asymptotic distribution of this estimator through a central limit theorem. 

The remainder if this article is organized as follows: in the next section, we recall the definition of multistable processes and our two examples of interest. We present the two estimators in Section \ref{Construct}.
Our main results on the convergence of the estimators are described in Section \ref{MR}. Subsection \ref{MRalpha} present the case of the index of stability $\alpha$. 
In subsection \ref{MRH}, we state the result giving the convergence of the estimator of the local self-similarity function $H$, with a central limit theorem in the case of the L\'evy multistable motion. 
In Section \ref{InterRes}, we give intermediate results which are used in the proofs of the main theorems. Section \ref{exa} contains applications of our results to two examples and real electrocardiographic data.
We give in Section \ref{Assum} a list of technical conditions on the kernel of multistable processes that involve the consistency of the estimators.
Finally we gather all the proofs of the statements of this article in Section \ref{proo}.

\section{Model}

Let us recall the definition of a {\it localisable process} \cite{Fal5,Fal6}: $Y =\{Y(t): t \in \bbbr\}$ 
is said to be localisable at $u$ if there exists an $H(u) \in \bbbr$ 
and a non-trivial limiting process $Y_{u}'$ such that
\begin{equation}
\lim_{r \to 0}\frac{Y(u+rt) -Y(u)}{r^{H(u)}} = Y_{u}'(t),
\label{locform1}
\end{equation}
where the convergence is in finite dimensional distributions. When the limit exits, $Y_{u}'=\{Y_{u}'(t): t\in \bbbr\}$ is termed the {\it local form} or tangent process of $Y$ at $u$. 

\medskip

{\bf Ferguson-Klass-LePage series representation}

\medskip
\noindent
We define now the multistable processes using the Ferguson-Klass-LePage series representation, that are defined as ``diagonals'' of random fields that we described below. In the sequel, $(E,{\cal E},m)$ will be a measure space, and $U$ an open interval of the real line $\bbbr$. 
We consider ${\cal F}_{\alpha}(E,{\cal E}, m)= \{ f: f \mbox{ is measurable and } \|f\|_{\alpha} < \infty\},$
where $\|\,\|_{\alpha}$ is the quasinorm (or norm if $1<\alpha \leq 2$) given by $\|f\|_{\alpha} = \left( \int_E |f(x)|^{\alpha}m(dx)\right)^{1/\alpha}.$ We will assume that $m$ is either a finite or a $\sigma$-finite measure, depending on the circumstances.

Let $\alpha$ be a $C^1$ function defined on $U$ and ranging in $[c,d] \subset (0,2)$. Let $f(t,u,.)$ be a family of functions such that, for all $(t,u) \in U^2$, $f(t,u,.) \in {\cal F}_{\alpha(u)}(E,{\cal E}, m)$.
We define also $r : E \rightarrow \mathbb{R}_+$ such that $\hat{m}(dx)=\frac{1}{r(x)}m(dx)$ is a probability measure. 
$(\Gamma_i)_{i \geq 1}$ will be a sequence of arrival times of a standard Poisson process and $(\gamma_i)_{i \geq 1}$  a sequence of i.i.d. random variables with distribution $P(\gamma_i=1)=P(\gamma_i=-1)=1/2$. Let $(V_i)_{i \geq 1}$ a sequence of i.i.d. random variables with distribution $\hat m$ on $E$ and we assume that the 
three sequences $(\Gamma_i)_{i \geq 1}$, $(V_i)_{i \geq 1}$, and $(\gamma_i)_{i \geq 1}$ are mutually independent. As in \cite{LGLV}, we will consider the following random field:

\begin{equation}\label{msfm2}
X(t,u)= C^{1/\alpha(u)}_{\alpha(u)} \sum_{i=1}^{\infty} \gamma_i \Gamma_i^{-1/\alpha(u)} r(V_i)^{1/\alpha(u)}f(t,u,V_i),
\end{equation}
where $C_{\eta} = \left( \int_{0}^{\infty} x^{-\eta} \sin (x)dx \right)^{-1}$.

Note that when the function $\alpha$  is constant, then \eqref{msfm2}
is just the Ferguson - Klass - LePage series representation of a stable
random variable (see \cite{BJP,FK,LP1,LP2,JR}
and \cite[Theorem 3.10.1]{ST} for specific properties of this representation).

\medskip

{\bf Multistable processes}

\medskip
\noindent
Multistable processes are obtained by taking diagonals on $X$ defined in (\ref{msfm2}), {\it i.e.}
\begin{equation}\label{Multdef}
Y(t) = X(t,t).
\end{equation}
 Indeed, as shown in Theorems 3.3 and 4.5 of \cite{LGLV}, provided some conditions are satisfied both by $X$ and by the function $f$, $Y$ will be
a localisable process whose local form is a stable process. We will always assume that $X(t,u)$ (as a process in $t$) is localisable at $u$ with exponent $H(u) \in (H_-,H_+) \subset (0,1)$, with local form $X'_u(t,u)$, and  $u \mapsto H(u)$ is a $C^1$ function.

We take as examples of multistable processes the ``multistable versions'' of some classical processes: the $\alpha$-stable L\'evy motion and the Linear Fractional Stable Motion.
In the sequel, $M$ will denote 
a symmetric $\alpha$-stable ($ 0 < \alpha < 2$)
random measure on $\bbbr$ with control measure Lebesgue measure ${\cal L}$.
We will write
$$L_{\alpha} (t) := \int_{0}^{t} M(dz)$$ 
for $\alpha$-stable L\'{e}vy motion, and we will use the Ferguson-Klass-LePage representation, 
\begin{displaymath}
 \forall t \in (0,1), \hspace{0.2cm} L_{\alpha} (t) = C_{\alpha}^{1/ \alpha} \sum_{i=1}^{\infty} \gamma_i \Gamma_i^{-1/\alpha}  \mathbf{1}_{[0,t]}(V_i).
\end{displaymath}
\noindent Let $\alpha: [0,1] \to (0,2)$ be continuously differentiable. Define 
\begin{displaymath}
 X(t,u)=C_{\alpha(u)}^{1/ \alpha(u)} \sum_{i=1}^{\infty} \gamma_i \Gamma_i^{-1/\alpha(u)}  \mathbf{1}_{[0,t]}(V_i)
 \end{displaymath}
and the symmetric multistable L\'{e}vy motion
\begin{displaymath}
Y(t) = X(t,t)=C_{\alpha(t)}^{1/ \alpha(t)} \sum_{i=1}^{\infty} \gamma_i \Gamma_i^{-1/\alpha(t)}  \mathbf{1}_{[0,t]}(V_i).
\end{displaymath}

The second example is a multistable version of the {\it well-balanced linear fractional $\alpha$-stable motion}:
\begin{displaymath}
L_{\alpha,H} (t) = \int_{-\infty}^{\infty} f_{\alpha,H}(t,x) M(dx) 
\end{displaymath}
where $t \in \mathbb{R}$, $H \in (0,1)$,and
\begin{displaymath}
f_{\alpha,H}(t,x) =  |t-x|^{H - 1/ \alpha} - |x|^{H - 1/ \alpha}.
\end{displaymath}

Let   $\alpha:\bbbr \to (0,2)$ and $H:\bbbr \to (0,1)$ be continuously differentiable. Define
\begin{equation} \label{fieldLmmm}
X(t,u)=C_{\alpha(u)}^{1/ \alpha(u)} \sum_{i,j=1}^{\infty} \gamma_i \Gamma_i^{-1/\alpha(u)} (|t-V_i|^{H(u)-1/ \alpha(u)}-|V_i|^{H(u)-1/\alpha(u)}) (\frac{\pi^2 j^2}{3})^{1/ \alpha(u)} \mathbf{1}_{[-j,-j+1[ \cup [j-1,j[}(V_i)
\end{equation}
and the linear multistable multifractional motion
\begin{equation}\label{defLmmm}
Y(t) = X(t,t).
\end{equation}
The localisability of L\'{e}vy motion and
linear fractional $\alpha$-stable motion simply stems from the fact that
they are self-similar with stationary increments \cite{Fal6}. We will apply our results to these processes,
that were defined in \cite{FLGLV,FLV}, in Section \ref{exa}.

\section{Construction of the estimators}\label{Construct}

Let $Y$ be a multistable process defined in (\ref{Multdef}). The estimation of the localisability function $H$ and the stability function $\alpha$ is based on the increments $(Y_{k,N})$ of $Y$. 
Define the sequence $(Y_{k,N})_{k \in \mathbb{Z}, N \in \mathbb{N}}$ by 
\begin{displaymath}
 Y_{k,N} = Y(\frac{k+1}{N})-Y(\frac{k}{N}).
\end{displaymath}

Let $t_0 \in \mathbb{R}$ be fixed. We introduce an estimator of $H(t_0)$ with 
\begin{displaymath}
 \hat{H}_N(t_0) = - \frac{1}{n(N)\log N}\sum_{k=[Nt_0]-\frac{n(N)}{2}}^{[Nt_0]+\frac{n(N)}{2}-1} \log |Y_{k,N}|
\end{displaymath}
where $(n(N))_{N \in \mathbb{N}}$ is a sequence taking even integer values. We expect the sequence $(\hat{H}_N(t_0))_N$ to converge to $H(t_0)$ thanks to the localisability of the process $Y$.
For the integers $k$ and $N$ such that $\frac{k}{N}$ is close to $t_0$, $\dfrac{Y_{k,N}}{(\frac{1}{N})^{H(t_0)}}$ is asymptotically distributed as $Y'_{t_0}(1)$.
More precisely $-\dfrac{\log |Y_{k,N} |}{\log N} = H(t_0) + \dfrac{Z_{k,N}}{\log N}$ where $(Z_{k,N})_{k,N}$ converge weakly to $-\log |Y'_{t_0}(1)|$ when $N$ tends to infinity and $\dfrac{k}{N}$ tends to $t_0$. 
We regulate the sequence $(Z_{k,N})$ near $t_0$ using the mean $\frac{1}{n(N)} \sum_{k=[Nt_0]-\frac{n(N)}{2}}^{[Nt_0]+\frac{n(N)}{2}-1}\limits Z_{k,N} $ and we can expect this sum will be bounded in the $L^r$ spaces to obtain the convergence with a rate $\frac{1}{\log N}$. The convergence is proved in Theorem \ref{ConvLpH}.

\bigskip

Let $p_0 \in (0,c)$ and $\gamma \in (0,1)$. With the increments of the process, we define the sample moments $S_N(p)$ by
\begin{displaymath}
 S_N(p)=\left( \frac{1}{n(N)} \sum_{k=[Nt_0]-\frac{n(N)}{2}}^{[Nt_0]+\frac{n(N)}{2}-1} |Y_{k,N}|^p \right)^{\frac{1}{p}}.
\end{displaymath}
Let \begin{displaymath}
 R_{\textrm{exp}}^{(N)}(p) = \frac{S_N(p_0)}{S_N(p)} \textrm{ and } R_{\alpha}(p) = \frac{(\E|Z|^{p_0})^{1/p_0}}{(\E|Z|^{p})^{1/p}} \mathbf{1}_{p < \alpha}
\end{displaymath}
where $Z$ is a standard symmetric $\alpha$-stable random variable (written $Z \sim S_{\alpha}(1,0,0)$ as in \cite{ST}), i.e $\E|Z|^{p} = \frac{2^{p-1}\Gamma(1-\frac{p}{\alpha})}{p\int_{0}^{+\infty} u^{-p-1} \sin^2 (u) du}$ . 
 
 Consider the set $A_N =: \argmin_{\alpha \in [0,2]}\limits \left( \int_{p_0}^{2} |R_{\textrm{exp}}^{(N)}(p) - R_{\alpha}(p) |^{\gamma} dp \right)^{1/ \gamma}.$ Since the function $\alpha \rightarrow \left( \int_{p_0}^{2}\limits |R_{\textrm{exp}}^{(N)}(p) - R_{\alpha}(p) |^{\gamma} dp \right)^{1/ \gamma}$ is a continuous function, $A_N$ is a non empty closed set. We define then an estimator of  $\alpha(t_0)$ by 
\begin{displaymath}
 \hat{\alpha}_N(t_0) = \min \left(\argmin_{\alpha \in [0,2]} \left( \int_{p_0}^{2} |R_{\textrm{exp}}^{(N)}(p) - R_{\alpha}(p) |^{\gamma} dp \right)^{1/ \gamma} \right).
\end{displaymath}
Under the conditions of Theorem \ref{ConvLpSnp}, $Y$ is $H(t_0)$-localisable and $Y'_{t_0}(1) \sim S_{\alpha(t_0)}(1,0,0)$ so $\dfrac{|Y_{k,N}|^p}{(\frac{1}{N})^{pH(t_0)}}$ converge weakly to $|Y'_{t_0}(1)|^p$ and with a meaning effect, $N^{H(t_0)} S_N(p)$ tends to $(\E|Y'_{t_0}(1)|^{p})^{1/p}$ in probability, which is the result of Theorem \ref{ConvLpSnp}. 
Following this, $ \int_{p_0}^{2} |R_{\textrm{exp}}^{(N)}(p) - R_{\alpha}(p) |^{\gamma} dp $ tends to $ \int_{p_0}^{2} |R_{\alpha(t_0)}(p) - R_{\alpha}(p) |^{\gamma} dp $ . 
Naturally, $\alpha(t_0)$ is the only solution of $\argmin_{\alpha \in [0,2]} \int_{p_0}^{2} |R_{\alpha(t_0)}(p) - R_{\alpha}(p) |^{\gamma} dp$ and this leads to the definition of $\hat{\alpha}_N(t_0) $. The consistency of $\hat{\alpha}_N(t_0) $ is proved in Theorem \ref{ConvLpalpha}.

\section{Main results}\label{MR}

The following theorems apply to a diagonal process $Y$ defined from the field $X$ given by (\ref{msfm2}). For convenience, the conditions required on $X$ and the function $f$ that appears in (\ref{msfm2}) are gathered in Section \ref{Assum}.
Theorem \ref{ConvLpalpha} leads to the convergence in the $L^r$ spaces of the estimator of the stability function $\alpha$, while Theorem \ref{ConvLpH} yield the convergence of the estimator of the localisability function $H$.
 We obtain also the convergence speed in the specific case of the symmetric multistable L\'evy motion. 

\subsection{Approximation of the stability function}\label{MRalpha}

\begin{theo}\label{ConvLpalpha}
  Let $Y$ be a multistable process and $t_0 \in U$. Assume the conditions (R1), (M1), (M2) and (M3). Assume in addition that:
 \begin{itemize}
 \item $\lim_{N \rightarrow +\infty}\limits n(N) = +\infty$ and $\lim_{N \rightarrow +\infty}\limits \frac{N}{n(N)} = +\infty$.
 \item The process $X(.,t_0)$ is $H(t_0)$-self-similar with stationary increments and $H(t_0) <1$.
\item $\lim_{j \rightarrow +\infty}\limits \int_E | h_{0,t_0}(x) h_{j,t_0}(x)|^{\frac{\alpha(t_0)}{2}} m(dx)=0$, where $h_{j,t_0}(x) = f(j+1,t_0,x)-f(j,t_0,x)$.
\end{itemize}
  Then for all $r > 0$,
 \begin{displaymath}
 \lim_{N \rightarrow +\infty}\E \left|\hat{\alpha}_N(t_0)- \alpha(t_0) \right|^r = 0.
\end{displaymath}
If, in addition, the conditions hold for all $t_0 \in U$, then for all $p >0$,

 \begin{displaymath}
 \lim_{N \rightarrow +\infty}\limits \E \left[ \int_U\limits |\hat{\alpha}_N(t)- \alpha(t) |^p dt\right] = 0.
\end{displaymath}

\end{theo}

{\bf Proof}
 
 See Section {\ref{proo}}.

\subsection{Approximation of the localisability function}\label{MRH}

\begin{theo}\label{ConvLpH}
  Let $Y$ be a multistable process. Assume that the localisability function $H$ and the function $\alpha$ are satisfying all the conditions (R1), (M1)-(M7) and (H1)-(H5) for an open interval $U$, and that $\lim_{N \rightarrow +\infty}\limits \frac{n(N)}{N} = 0.$\par
   Then, for all $t_0 \in U$ and all $r > 0$,
 \begin{displaymath}
 \lim_{N \rightarrow +\infty}\E \left|\hat{H}_N(t_0)- H(t_0) \right|^r = 0.
\end{displaymath}

Moreover, for all $[a,b] \subset U$ and all $p >0$,

 \begin{displaymath}
 \lim_{N \rightarrow +\infty}\limits \E \left[ \int_a^b\limits |\hat{H}_N(t)- H(t) |^p dt\right] = 0.
\end{displaymath}

\end{theo}

{\bf Proof}
 
 See Section {\ref{proo}}.

\medskip

{\bf Remark:} Under the conditions (R1), (M1), (M2) and (M3) listed in the theorem, Theorems 3.3 and 
4.5 of \cite{LGLV} imply that $Y$ is $H(t_0)-$localisable at $t_0$.

 We obtain for the symmetric multistable L\'{e}vy motion the convergence in distribution of the estimator $\hat{H}_N(t_0)$ in the following theorem. 
We expect the same result holds for a more general class of processes, in particular when the conditions of Theorem \ref{ConvReste} are satisfied. For $Z$ a standard $\alpha(t_0)$-stable random variable, we define $\mu_{t_0} = \E [\log |Z|]$ and $\sigma^2_{t_0}= Var(\log |Z| )$.

 \begin{theo}\label{Convlaw}
  Let $Y$ be a symmetric multistable L\'evy motion with $\alpha : [0,1] \rightarrow (1,2)$ continuously differentiable, and $t_0 \in (0,1)$. Assume that $n(N) = O( N^{\delta})$ with $\delta \in \left(0, \frac{2\alpha(t_0)-2}{3\alpha(t_0)+2} \right)$. Then
  
  \begin{displaymath}
   \sqrt{n(N)} \left( \log N \left( \hat{H}_N(t_0)- H(t_0) \right) + \mu_{t_0} \right) \cd \mathcal{N}(0, \sigma^2_{t_0} )
  \end{displaymath}
as $N \rightarrow +\infty$.
 \end{theo}

\section{Intermediate results}\label{InterRes}

All the proofs of the intermediate results are stated in Section \ref{proo}.
We first give conditions for the convergence in probability of $S_N(p)$ in Theorem \ref{ConvLpSnp}, which is useful to establish the consistency of the estimator $\hat{\alpha}_N(t_0)$.

\begin{theo}\label{ConvLpSnp}
  Let $Y$ be a multistable process. Assume the conditions (R1), (M1), (M2) and (M3). Assume in addition that:
 \begin{itemize}
 \item $\lim_{N \rightarrow +\infty}\limits n(N) = +\infty$ and $\lim_{N \rightarrow +\infty}\limits \frac{N}{n(N)} = +\infty$.
 \item The process $X(.,t_0)$ is $H(t_0)$-self-similar with stationary increments and $H(t_0) <1$.
\item $\lim_{j \rightarrow +\infty}\limits \int_E | h_{0,t_0}(x) h_{j,t_0}(x)|^{\frac{\alpha(t_0)}{2}} m(dx)=0$, where $h_{j,t_0}(x) = f(j+1,t_0,x)-f(j,t_0,x)$.
\end{itemize}

 Then, for all $p \in [p_0,\alpha(t_0))$,
 \begin{displaymath}
  N^{H(t_0)}S_N(p) \underset{N \rightarrow +\infty}{\longrightarrow} (\E|X(1,t_0)|^{p})^{1/p}
\end{displaymath}
where the convergence is in probability.
\end{theo}

\medskip

\noindent We establish under several assumptions that the sequence $(\hat{H}_N(t))_{N}$ is almost surely uniformly bounded on every compact $[a,b] \subset U$. 

\begin{lem}\label{LemSup}
 Assume that the localisability function $H$ and the function $\alpha$ are satisfying all the conditions (R1), (M1)-(M7) and (H1)-(H5) for an open interval $U$, and that $\lim_{N \rightarrow +\infty}\limits \frac{n(N)}{N} = 0.$ Then there exists $B \in \bbbr$ such that for all $[a,b] \subset U$, 
 \begin{displaymath}
  \P \left( \liminf_{N \rightarrow + \infty}\limits \{ \sup_{t \in [a,b]}\limits | \hat{H}_N(t)| \leq B\} \right) = 1.
 \end{displaymath}

\end{lem}

 \medskip
 
\noindent We state then a theorem implying the rate of convergence of the estimator $\hat{H}_N(t_0)$.
 
 \begin{theo}\label{ConvReste}
  Let $Y$ be a multistable process and $t_0 \in U$. Assume the conditions (R1), (M1), (M2), (M3). Assume in addition that :
  \begin{itemize}
   \item $n(N) = O( N^{\delta})$ with $\delta \in \left(0, \frac{2\alpha(t_0)(1-H(t_0))}{2+3\alpha(t_0)} \right)$,
   \item The process $X(.,u)$ is $H(u)$-self-similar with stationary increments and $H(u) <1$, for all $u \in U$.
  \end{itemize}
Then
\begin{displaymath}
 \lim_{N \rightarrow +\infty}\limits \frac{1}{\sqrt{n(N)}} \sum_{k=[Nt_0]-\frac{n(N)}{2}}^{[Nt_0] + \frac{n(N)}{2}-1}\limits \log \left| \frac{Y(\frac{k+1}{N}) -Y(\frac{k}{N})}{X(\frac{k+1}{N},\frac{k}{N}) - X(\frac{k}{N},\frac{k}{N})}\right| =0
\end{displaymath}
where the convergence is in probability.
 \end{theo}

 \medskip
 
 \noindent Finally, we set up a technical lemma, which is useful for Theorem \ref{ConvReste}. 
 
 \begin{lem}\label{ProbaRapport}
  Assume the conditions (R1), (M1), (M2) and (M3). Let $t_0 \in U$. If $X(.,u)$ is $H(u)$-self-similar with stationary increments and $H(u) <1$, for all $u \in U$, then there exists $K_U >0$ such that for all $\lambda \in (0,1/e)$, for all $(k,N) \in \bbbz \times \bbbn$ with $k \in \left[ [Nt_0]-\frac{n(N)}{2} , [Nt_0]+\frac{n(N)}{2} -1\right]$, 
\begin{displaymath}
\P \left( \frac{|X(\frac{k+1}{N},\frac{k+1}{N}) - X(\frac{k+1}{N},\frac{k}{N}) |}{|X(\frac{k+1}{N},\frac{k}{N}) - X(\frac{k}{N},\frac{k}{N}) |} > \lambda \right) \leq K_U  \frac{|\log N|^d|\log \lambda|^d}{N^{\frac{d(1-H_{-})}{1+c}} \lambda^{\frac{d}{1+c}}}.
\end{displaymath}
 \end{lem}

\section{Examples and simulations}\label{exa}

In this section, we apply the results to our two examples: the Linear multifractional multistable motion and the multistable L\'evy motion. 
We provide then an example of application with ECG data.

\subsection{Linear multistable multifractional motion}

We consider first the Linear multistable multifractional motion (Lmmm) defined by (\ref{defLmmm}).

\begin{prop}
 Assume that $H- \frac{1}{\alpha}$ is a non-negative function, $\lim_{N \rightarrow +\infty}\limits n(N) = +\infty$ and $\lim_{N \rightarrow +\infty}\limits \frac{N}{n(N)} = +\infty$. Then for all $r >0$ and all $[a,b] \subset \bbbr$,

 \begin{displaymath}
  \lim_{N \rightarrow +\infty}\E \left[ \int_a^b\limits |\hat{\alpha}_N(t)- \alpha(t) |^r dt \right] = 0,
 \end{displaymath}
and for all $t_0 \in \bbbr$, 
\begin{displaymath}
  \lim_{N \rightarrow +\infty}\E \left|\hat{H}_N(t_0)- H(t_0) \right|^r = 0.
\end{displaymath} 
 \end{prop}

 {\bf Proof}
 
 Let $t_0 \in [a,b] \subset \bbbr$ and $r >0$.
 
  \noindent We know from \cite{LGLV3} that the conditions (R1), (M1), (M2) and (M3) are satisfied. 
 Since the process $X(.,t_0)$ is a $\left(H(t_0), \alpha(t_0) \right)$ linear fractional stable motion, $X(.,t_0)$ is $H(t_0)$-self-similar with stationary increments \cite{ST}.
 Let us show that  $\lim_{j \rightarrow +\infty}\limits \int_{\bbbr}\limits | h_{0,t_0}(x) h_{j,t_0}(x)|^{\frac{\alpha(t_0)}{2}}dx=0$.

 Let $\varepsilon > 0$. Let $c_0 >0$ such that $\int_{|x| > c_0}\limits  | h_{0,t_0} (x)|^{\alpha(t_0)} dx \leq \frac{\varepsilon}{2} $. By the Cauchy-Schwartz inequality, we have that 
 $$\int_{|x| > c_0}\limits  | h_{0,t_0} (x)h_{j,t_0}(x)|^{\alpha(t_0)/2} dx \leq (\frac{\varepsilon}{2})^{1/2} \| h_{j,t_0}\|_{\alpha(t_0)}^{\alpha(t_0)/2}=(\frac{\varepsilon}{2})^{1/2} \| h_{0,t_0}\|_{\alpha(t_0)}^{\alpha(t_0)/2}.$$
 This implies the desired convergence since $\forall x \in [-c_0,c_0]$,
 $\lim_{j \rightarrow +\infty}\limits | h_{0,t_0}(x) h_{j,t_0}(x)|^{\frac{\alpha(t_0)}{2}} =0$, $(h_{j,t_0}(x))_j$ is uniformly bounded on $[-c_0,c_0]$, and therefore
 \begin{displaymath}
  \lim_{j \rightarrow +\infty} \int_{|x| \leq c_0} | h_{0,t_0}(x) h_{j,t_0}(x)|^{\frac{\alpha(t_0)}{2}} dx=0.
 \end{displaymath}
  We deduce from Theorem \ref{ConvLpalpha} that for all $t_0 \in [a,b]$, $ \lim_{N \rightarrow +\infty}\limits \E \left|\hat{\alpha}_N(t_0)- \alpha(t_0) \right|^r = 0$.
 Since $\hat{\alpha}$ and $\alpha$ are bounded by $2$, $  \lim_{N \rightarrow +\infty}\E \left[ \int_a^b\limits |\hat{\alpha}_N(t)- \alpha(t) |^r dt \right] = 0$.
 
 Let $t_0 \in \bbbr$.
 We know from \cite{LGLV3} that there exists $U$ an open interval such that $t_0 \in U$ and (M3), (M4), (M5), (M6), (M7), (H1)-(H5) hold. 
  We deduce from Theorem \ref{ConvLpH} that
  \begin{displaymath}
     \lim_{N \rightarrow +\infty}\E \left|\hat{H}_N(t_0)- H(t_0) \right|^r = 0 \quad \Box
  \end{displaymath}
 
\begin{figure}[H]
\begin{tabular}{lll}
  \includegraphics[scale=0.28]{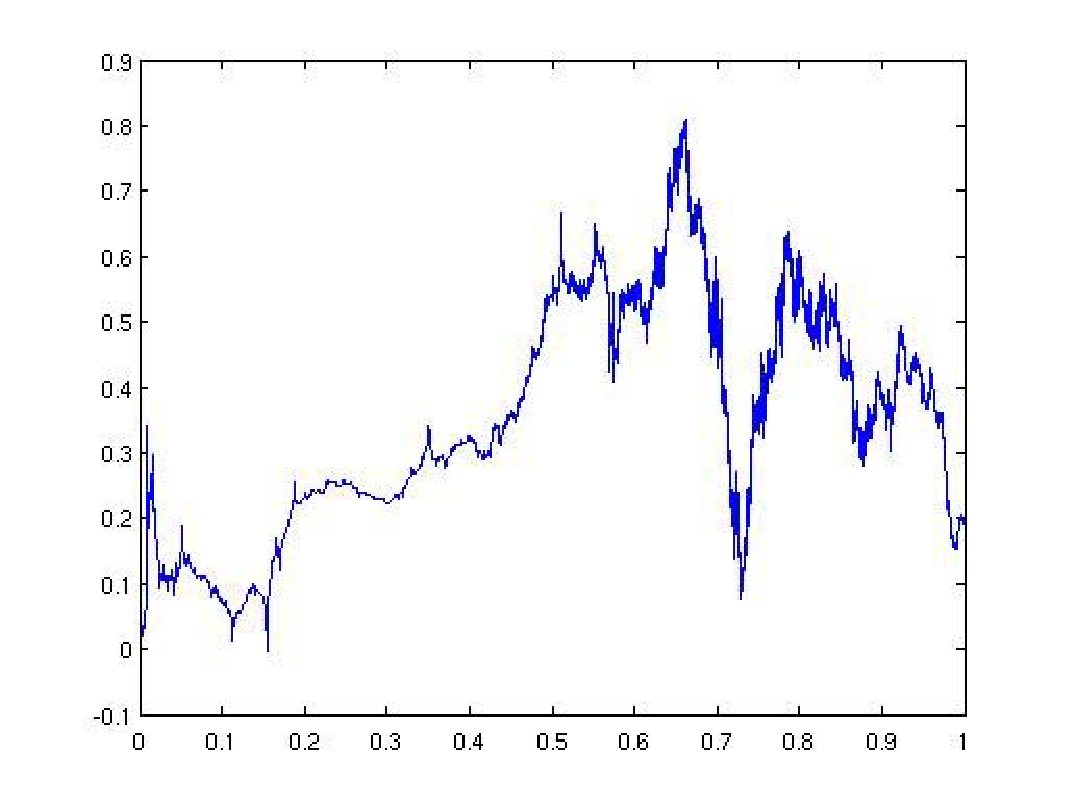} & \includegraphics[scale=0.28]{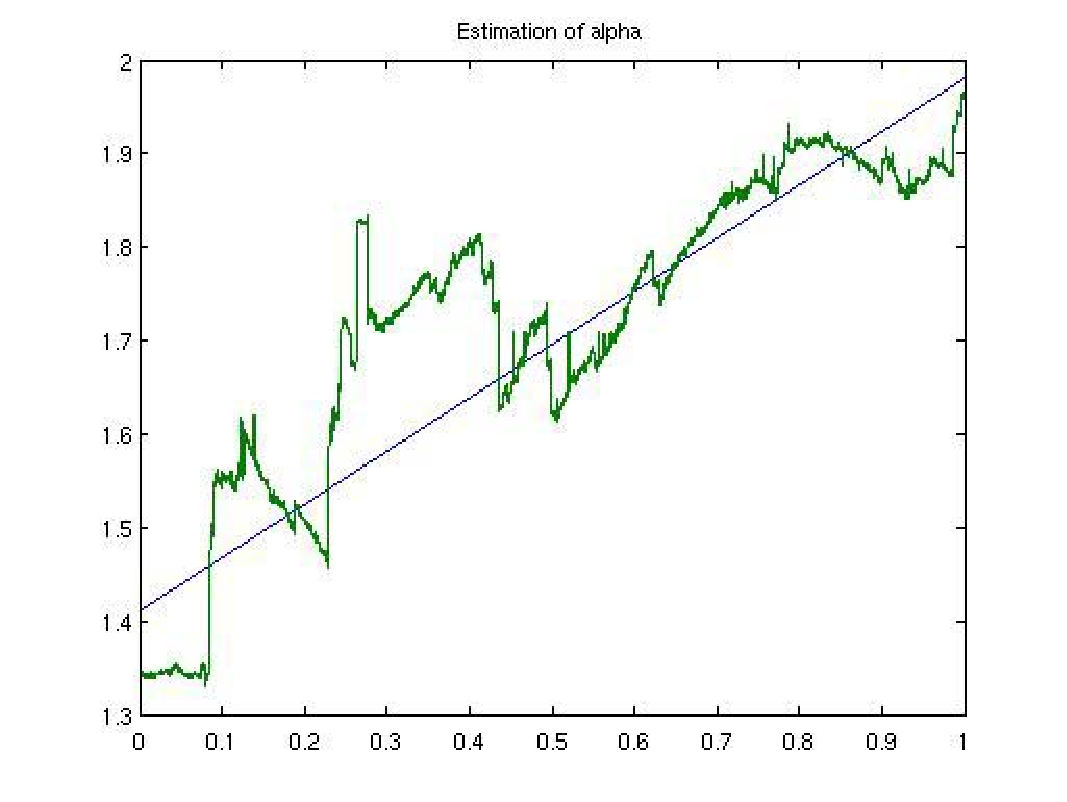} &
  \includegraphics[scale=0.28]{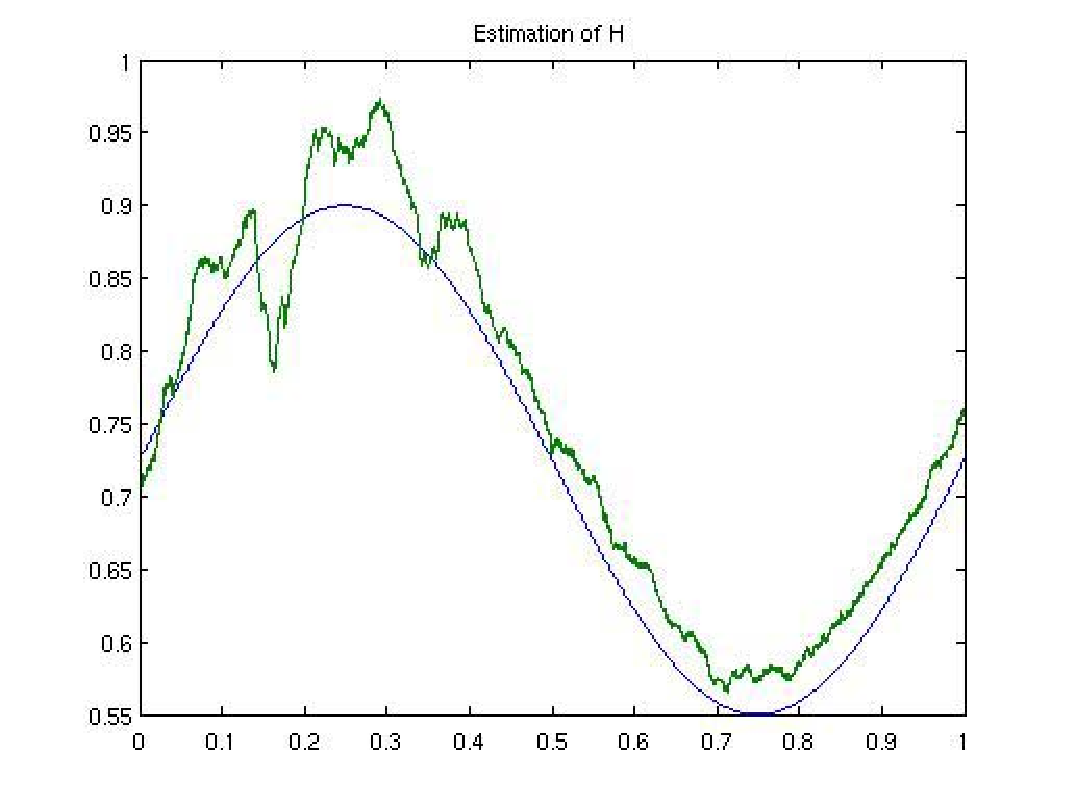}\\
& {\scriptsize  \hspace{1.5cm }$\alpha(t)=1.41+0.57t$} & {\scriptsize  \hspace{0.7cm } $H(t)=0.725+0.175\sin(2 \pi t)$ }\\
\includegraphics[scale=0.28]{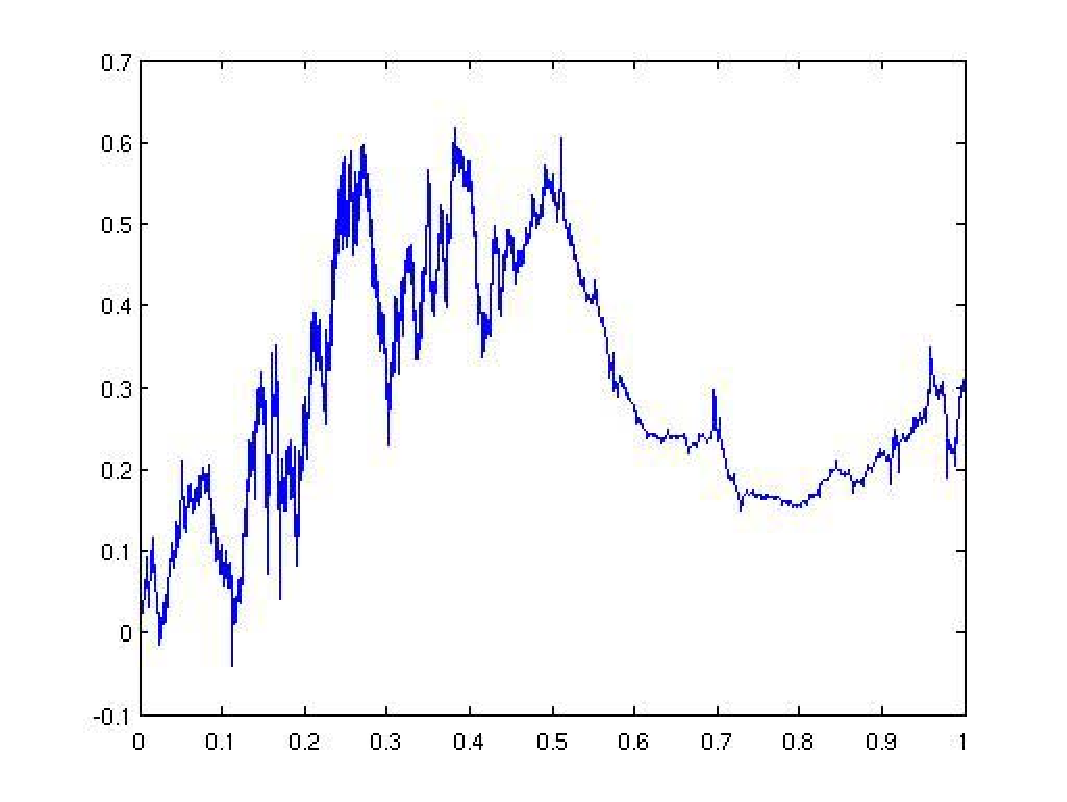} & \includegraphics[scale=0.28]{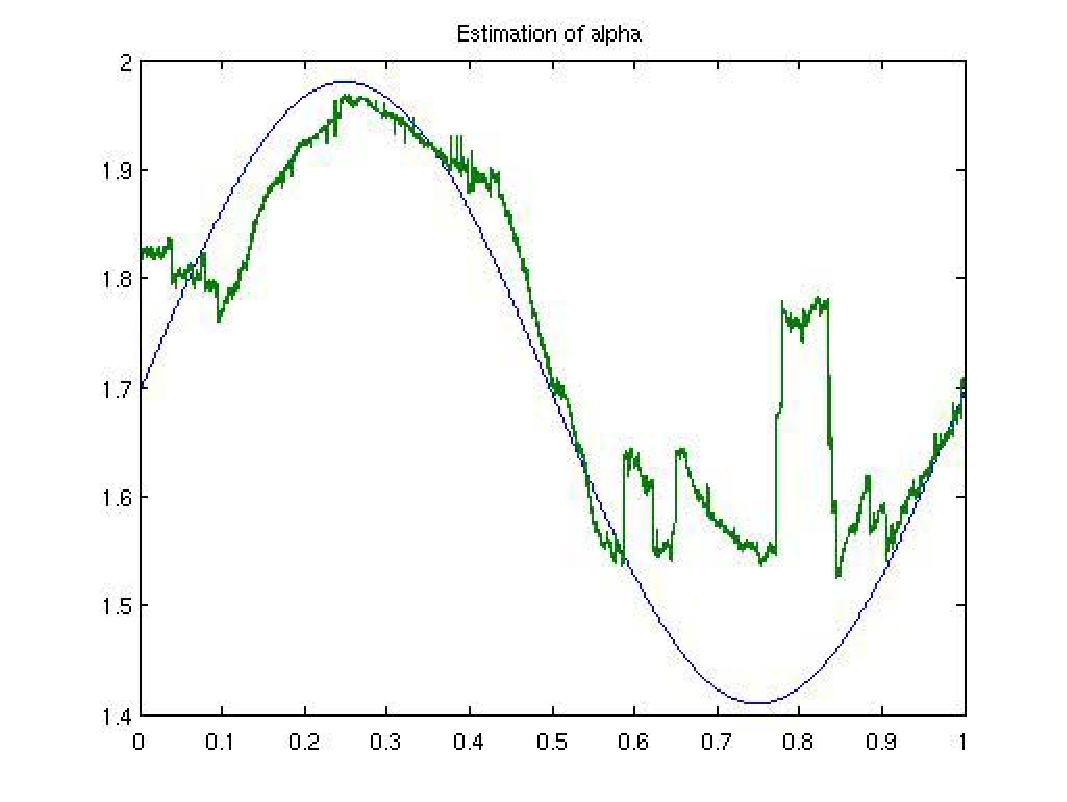} &
  \includegraphics[scale=0.28]{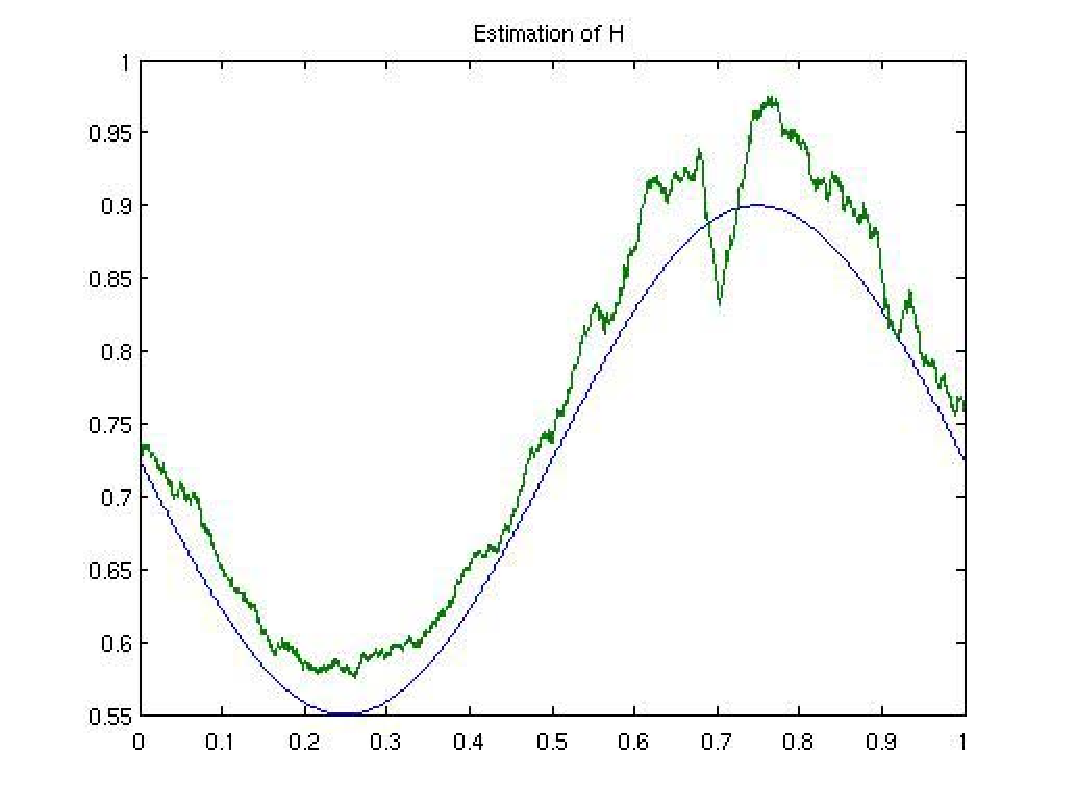}\\
& {\scriptsize  \hspace{0.8cm }$\alpha(t)=1.695+0.235 \sin(2\pi t)$} & {\scriptsize  \hspace{0.7cm } $H(t)=0.725-0.175\sin(2 \pi t)$ }\\
\includegraphics[scale=0.28]{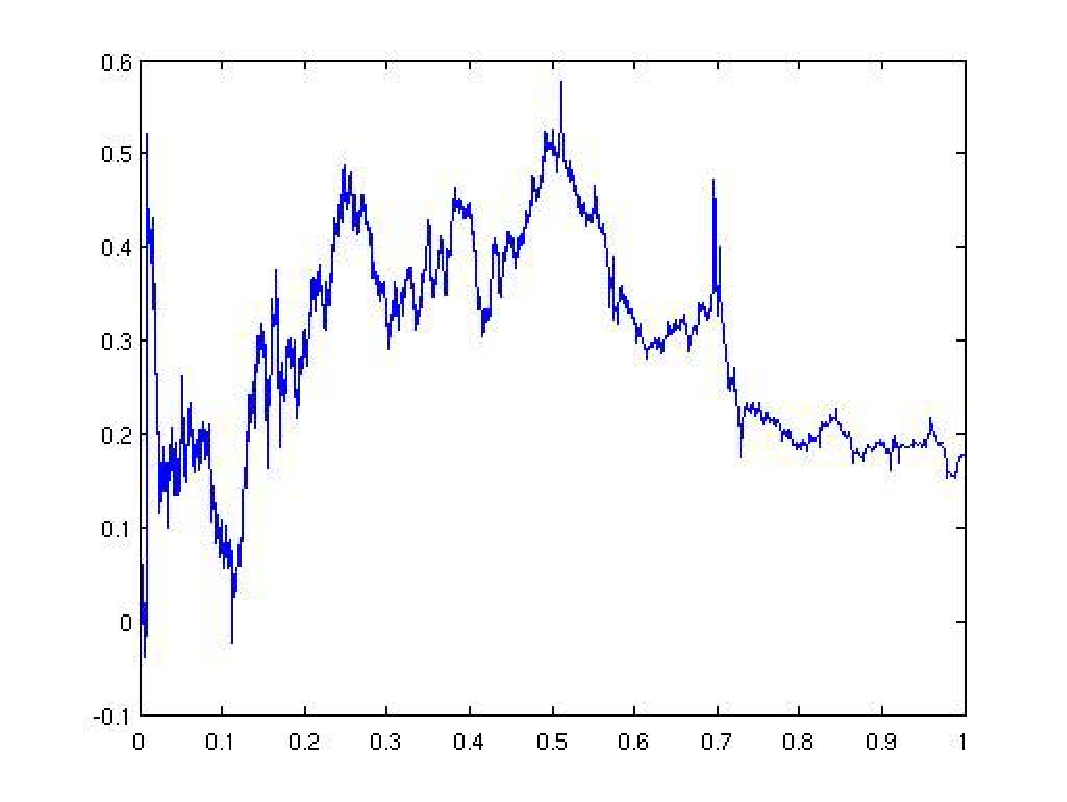} & \includegraphics[scale=0.28]{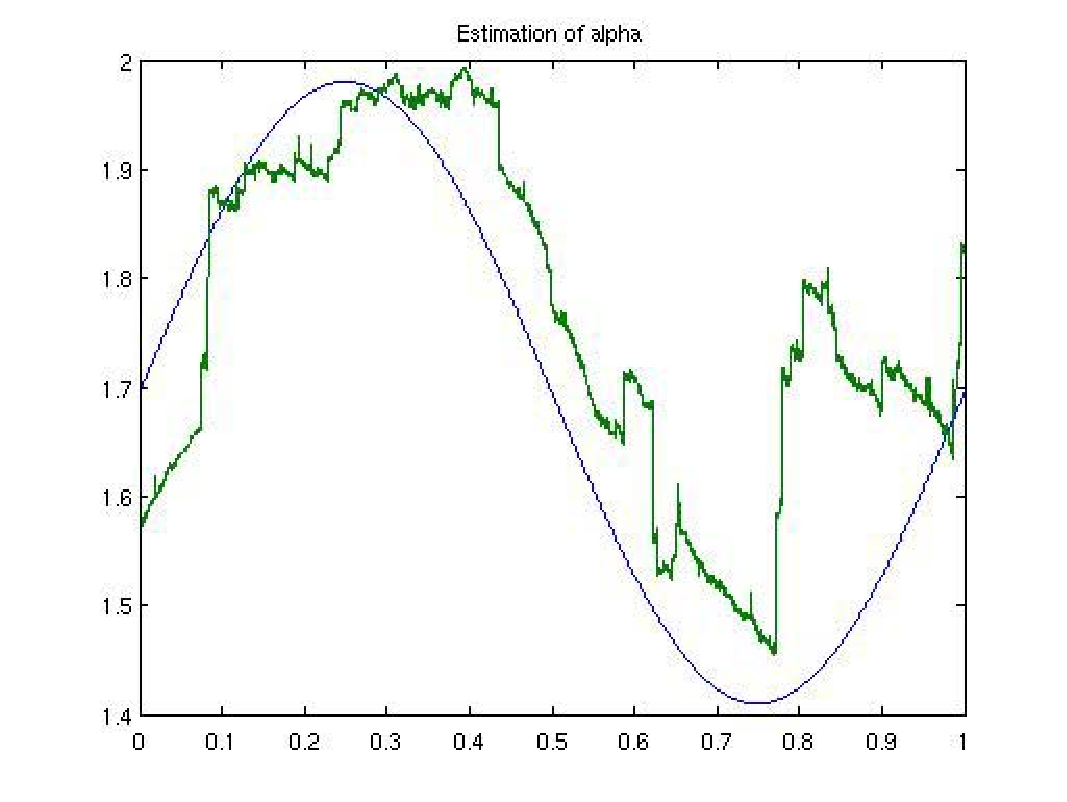} &
  \includegraphics[scale=0.28]{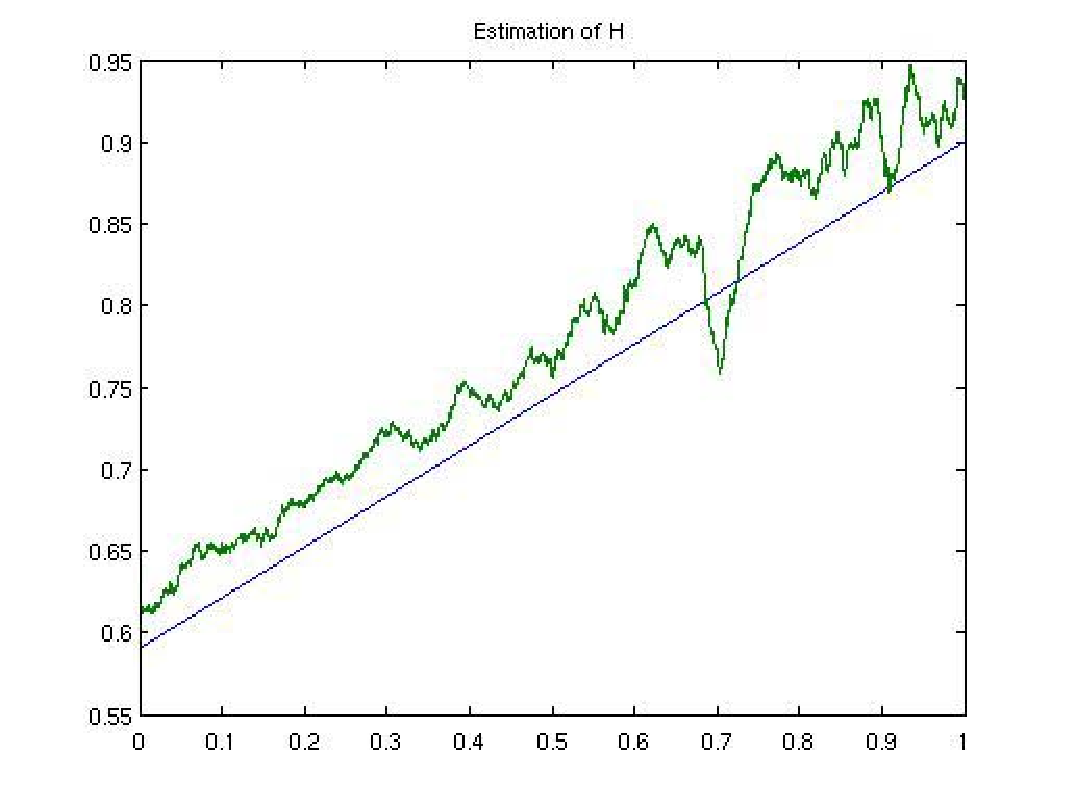}\\
& {\scriptsize  \hspace{0.8cm }$\alpha(t)=1.695+0.235 \sin(2\pi t)$} & {\scriptsize  \hspace{1.5cm } $H(t)=0.59+0.31 t$ }\\
\includegraphics[scale=0.28]{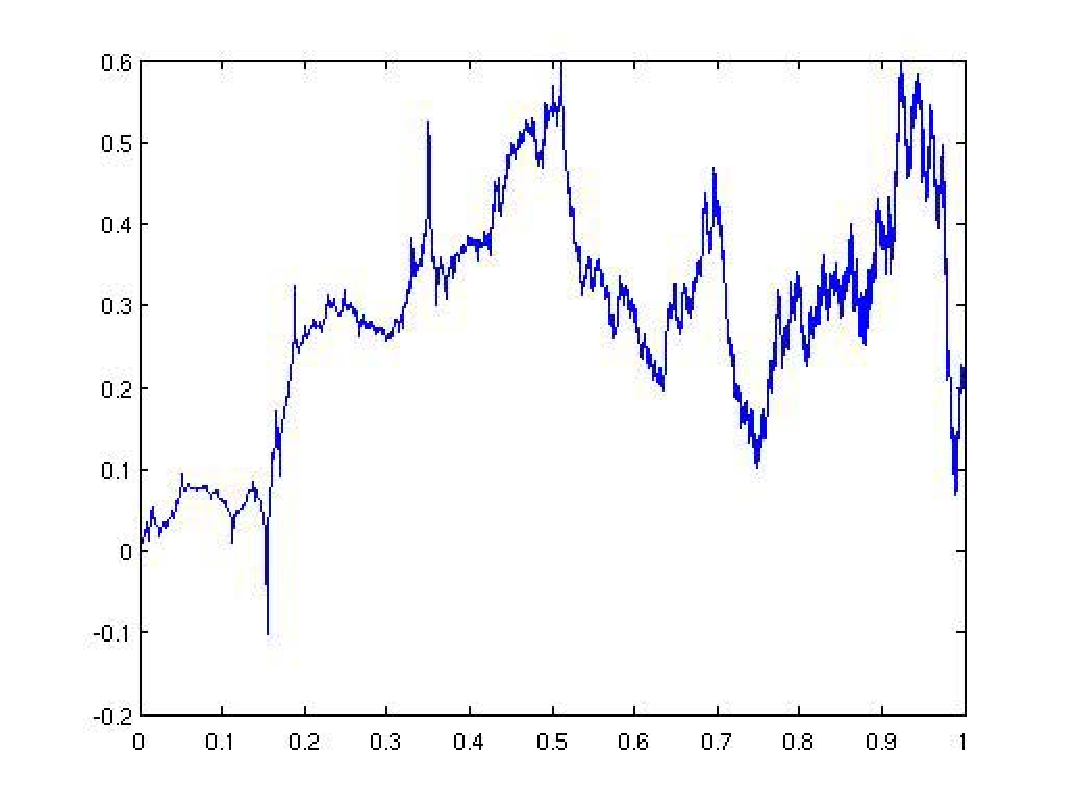} & \includegraphics[scale=0.28]{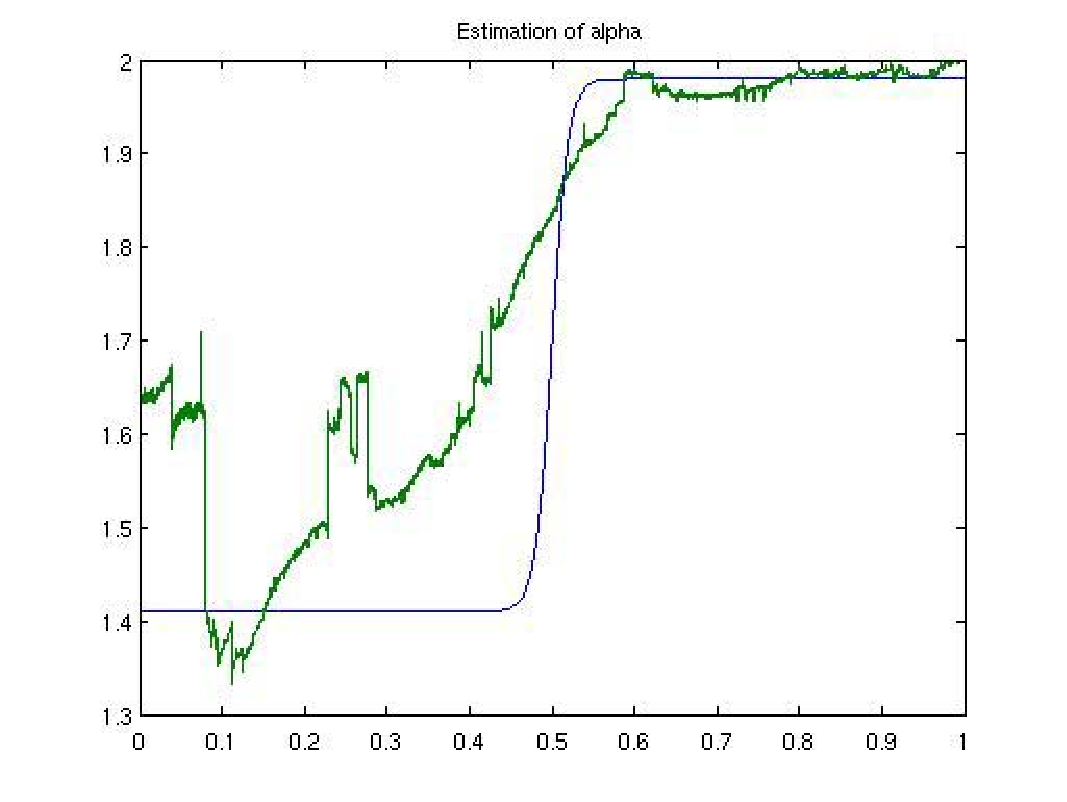} &
  \includegraphics[scale=0.28]{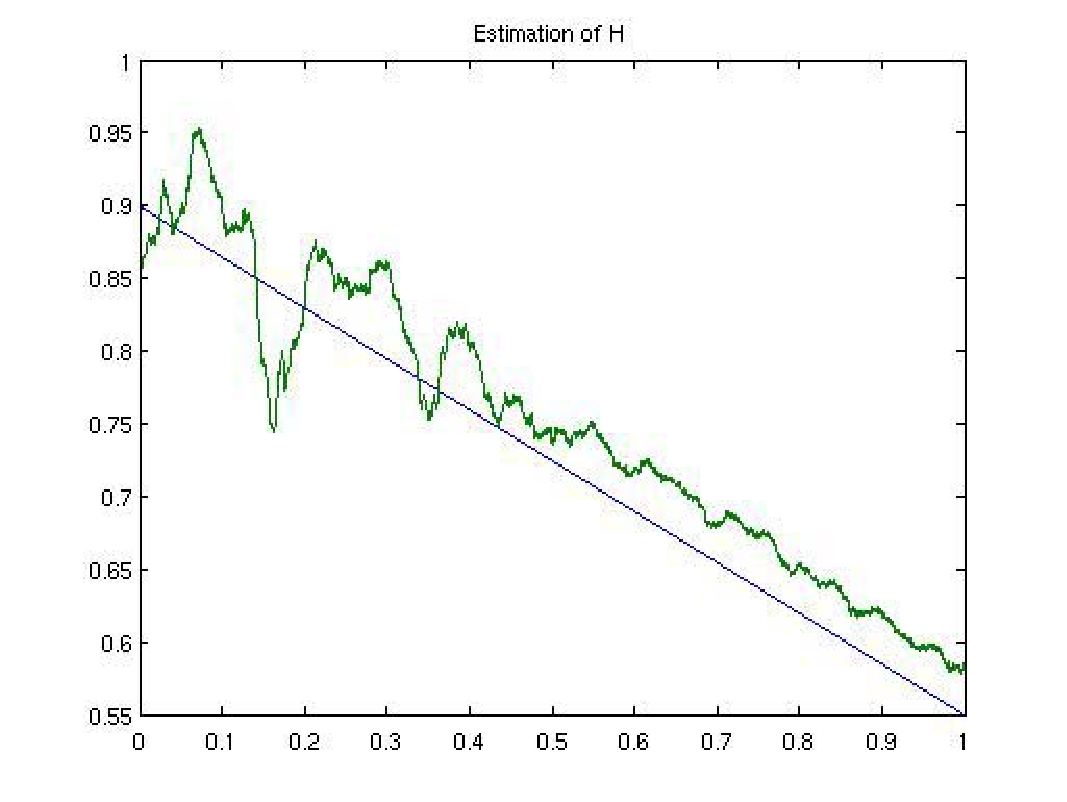}\\
& {\scriptsize  \hspace{0.6cm }$\alpha(t)=1.41+ \dfrac{0.47}{1+\exp(20-40 t)}$} & {\scriptsize  \hspace{1.4cm } $H(t)=0.9 - 0.35 t$ }\\
\includegraphics[scale=0.28]{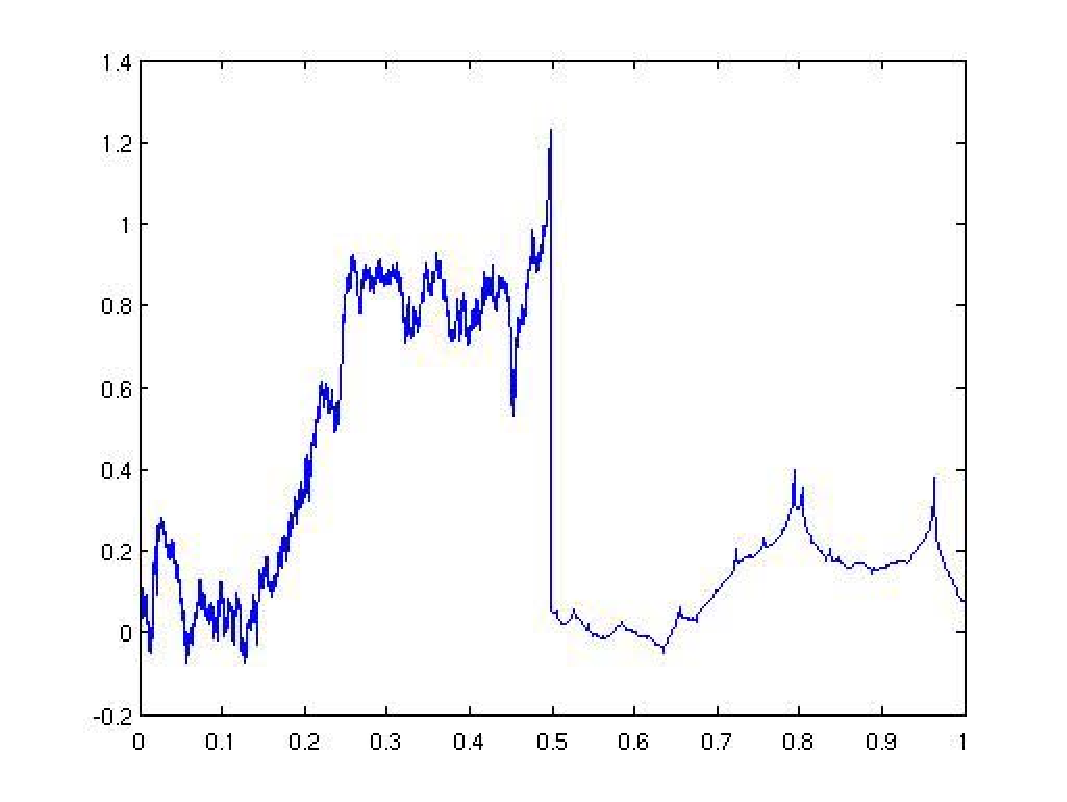} & \includegraphics[scale=0.28]{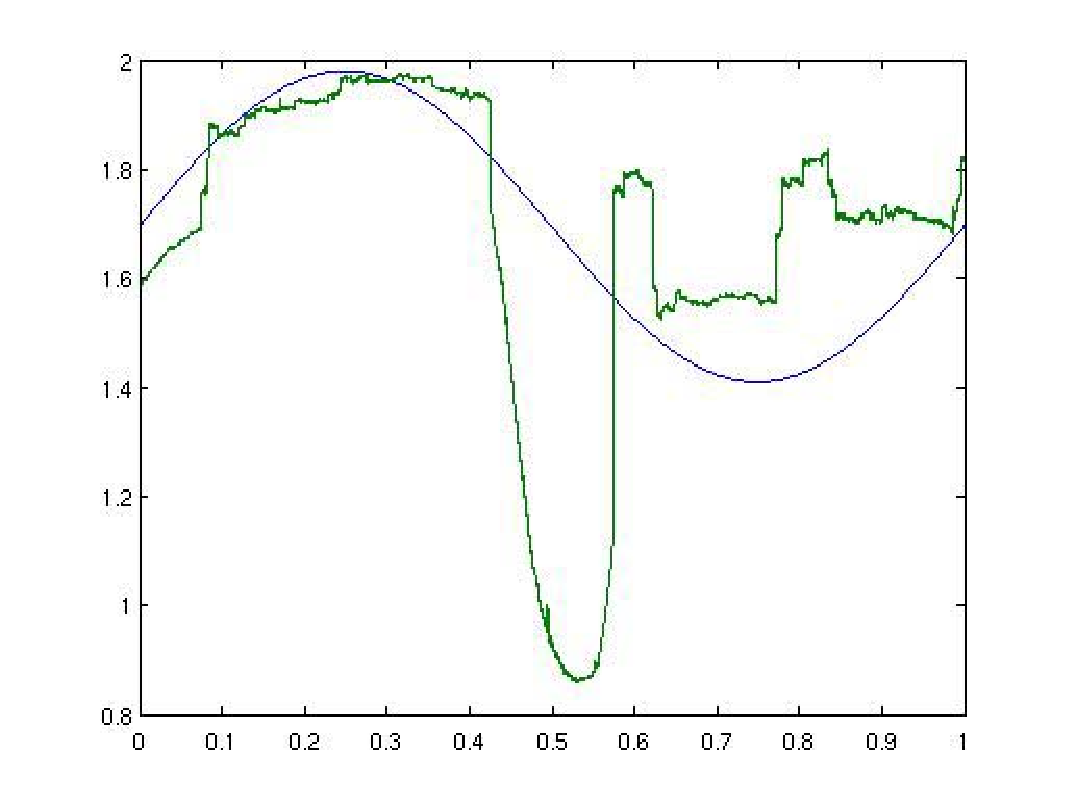} &
  \includegraphics[scale=0.28]{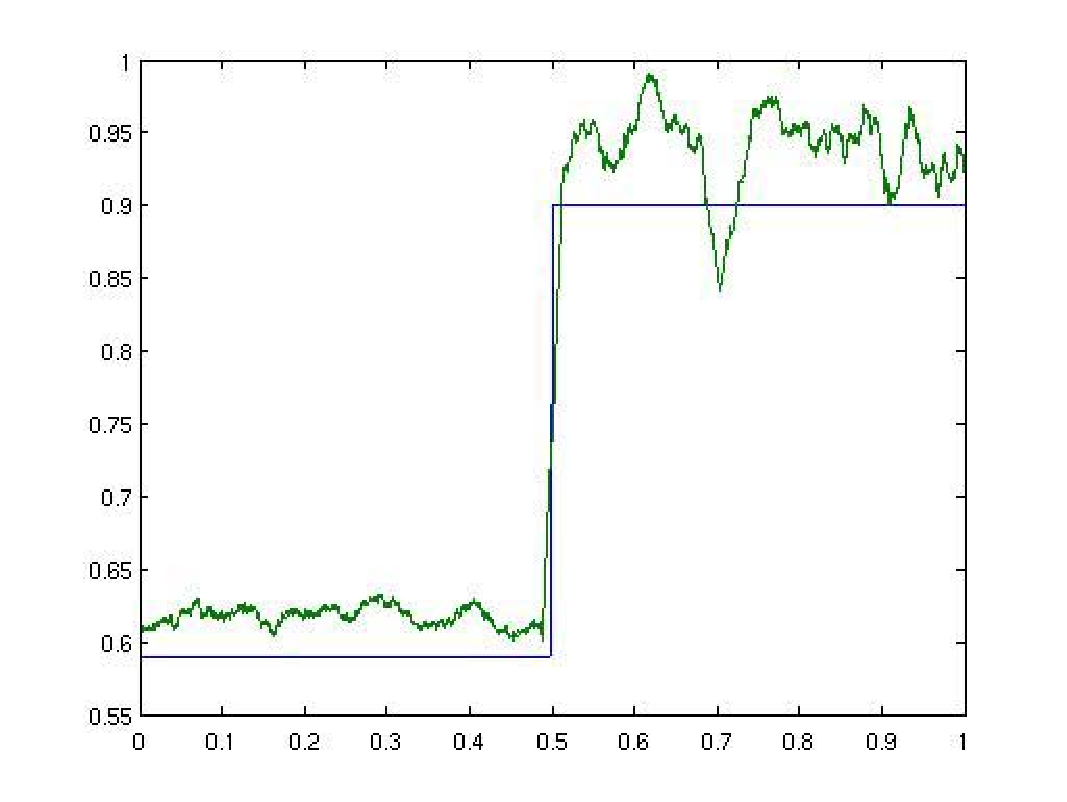}\\
  {\footnotesize \hspace{2.2cm}  a)} & {\footnotesize \hspace{2.2cm}  b)} & {\footnotesize \hspace{2.2cm}  c)}\\
\end{tabular}
\caption{Trajectories with $N=20000$ in the first column, the estimations of $\alpha$ with $n(N)=3000$ points in the second column, and in the last one, the estimations of $H$ with $n(N)=500$ points.}\label{fig4}
\end{figure}

  We show on Figure \ref{fig4} some paths of Lmmm, with the two corresponding estimations of $\alpha$ and $H$.
  To simulate the trajectories, we have used the field (\ref{fieldLmmm}). All the increments of $X(.,u)$ are $\left(H(u), \alpha(u) \right)$-linear fractional stable motions, generated using the LFSN program of \cite{ST3}. After we have taken the diagonal process $X(t,t)$.

These estimates are overall further than the estimates in the case of the Levy process, because of greater correlations between the increments of the process. 
However, the estimation of $H$ does not seem to be disturbed by this dependance.
The shape of the function $H$ is kept. For $\alpha$, we notice some disruptions when the function is close to $1$. We finally show an example where the estimation of $\alpha$ is not good enough in the last line of Figure \ref{fig4}.
The trajectory, Figure \ref{fig4}.a), seems to have a big jump, which leads to decrease the estimator $\hat{\alpha}$, represented on Figure \ref{fig4}.b), while the jump is taken account in the $n(N)$ points. 
The estimation of $H$, represented on Figure \ref{fig4}.c), does not seem to be affected by this phenomenon.
 

 \subsection{Symmetric multistable L\'{e}vy motion}

\noindent Let $\alpha: [0,1] \to (1,2)$ be continuously differentiable. Define 
\begin{equation}\label{fieldLevy}
 X(t,u)=C_{\alpha(u)}^{1/ \alpha(u)} \sum_{i=1}^{\infty} \gamma_i \Gamma_i^{-1/\alpha(u)}  \mathbf{1}_{[0,t]}(V_i)
 \end{equation}
and the symmetric multistable L\'{e}vy motion
\begin{displaymath}
Y(t) = X(t,t)=C_{\alpha(t)}^{1/ \alpha(t)} \sum_{i=1}^{\infty} \gamma_i \Gamma_i^{-1/\alpha(t)}  \mathbf{1}_{[0,t]}(V_i).
\end{displaymath}

\begin{prop}\label{ExaLevy}
If $\lim_{N \rightarrow +\infty}\limits n(N) = +\infty$ and $\lim_{N \rightarrow +\infty}\limits \frac{N}{n(N)} = +\infty$, then for all $r >0$,
 \begin{displaymath}
 \lim_{N \rightarrow +\infty}\E \left[ \int_0^1 |\hat{\alpha}_N(t)- \alpha(t)|^r dt \right] = 0.
\end{displaymath}
For all $[a,b] \subset (0,1)$, 
\begin{displaymath}
 \lim_{N \rightarrow +\infty}\E \left[ \int_a^b | \hat{H}_N(t)- \frac{1}{\alpha(t)} |^r dt \right] = 0.
\end{displaymath}
Let $t_0 \in (0,1)$. If we assume in addition that $n(N) = O( N^{\delta})$ with $\delta \in \left(0, \frac{2\alpha(t_0)-2}{3\alpha(t_0)+2} \right)$. Then
  
  \begin{displaymath}
   \sqrt{n(N)} \left( \log N \left( \hat{H}_N(t_0)- H(t_0) \right) + \mu_{t_0} \right) \cd \mathcal{N}(0, \sigma^2_{t_0} )
  \end{displaymath}
as $N \rightarrow +\infty$.

\end{prop}

{\bf Proof}

We know from \cite{LGLV3} that the conditions (R1), (M1), (M2) and (M3) are satisfied with $U=(0,1)$. 
 Since the process $X(.,t_0)$ is a L\'evy motion $\alpha(t_0)$-stable, $X(.,t_0)$ is $\frac{1}{\alpha(t_0)}$-self-similar with stationary increments \cite{ST}.
 $h_{j,t_0}(x) = \mathbf{1}_{[j,j+1[}(x)$ so for $j \geq 1$,
\begin{displaymath}
\int_{\bbbr} | h_{0,t_0}(x) h_{j,t_0}(x)|^{\frac{\alpha(t_0)}{2}} dx =0. 
\end{displaymath}
We conclude with Theorem \ref{ConvLpalpha} that $\lim_{N \rightarrow +\infty}\limits \E \left[ \int_0^1 |\hat{\alpha}_N(t)- \alpha(t)|^r dt \right] = 0. $.

Let $[a,b] \subset (0,1)$. 
We easily check that the nine conditions (M4)-(M7) and (H1)-(H5) are satisfied with $U=(a,b)$ and $H(t) = \frac{1}{\alpha(t)}$. We conclude with Theorem \ref{ConvLpH} that $ \lim_{N \rightarrow +\infty}\limits \E \left[ \int_a^b |\hat{H}_N(t)- \frac{1}{\alpha(t)}|^r dt \right] = 0.$
The end of Proposition \ref{ExaLevy} is a reminder of Theorem \ref{Convlaw} \Box

\begin{figure}[H]
\begin{tabular}{lll}
  \includegraphics[scale=0.3]{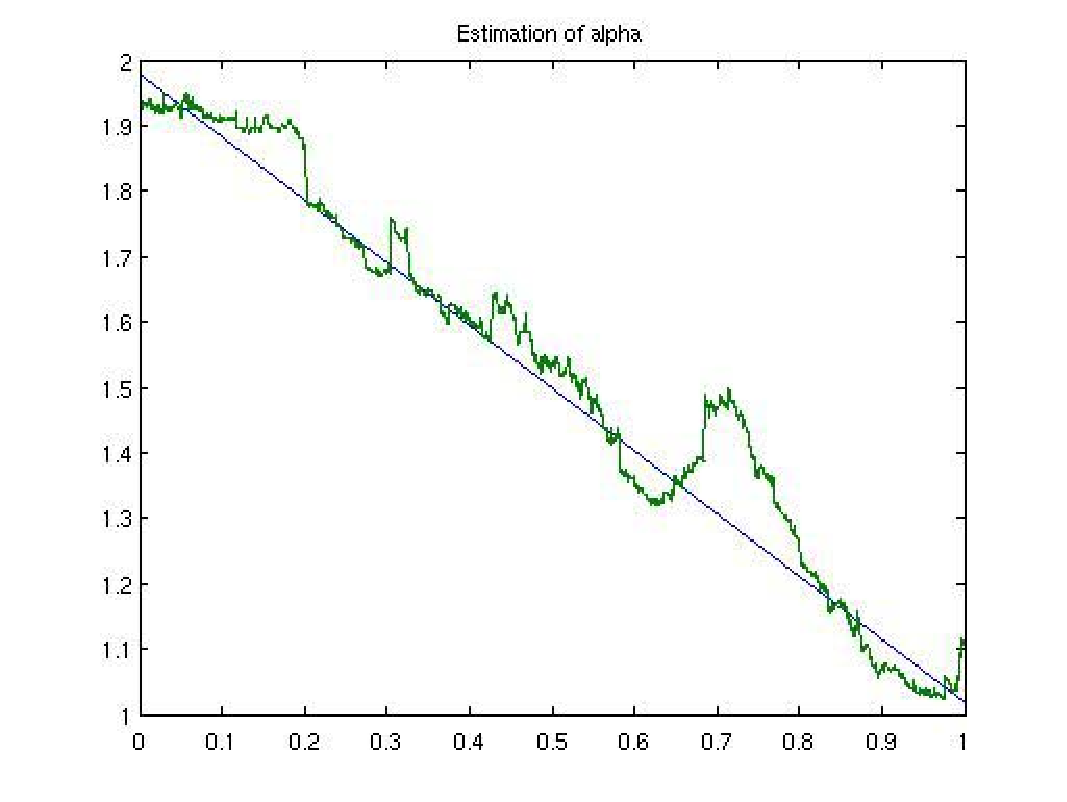} & \includegraphics[scale=0.3]{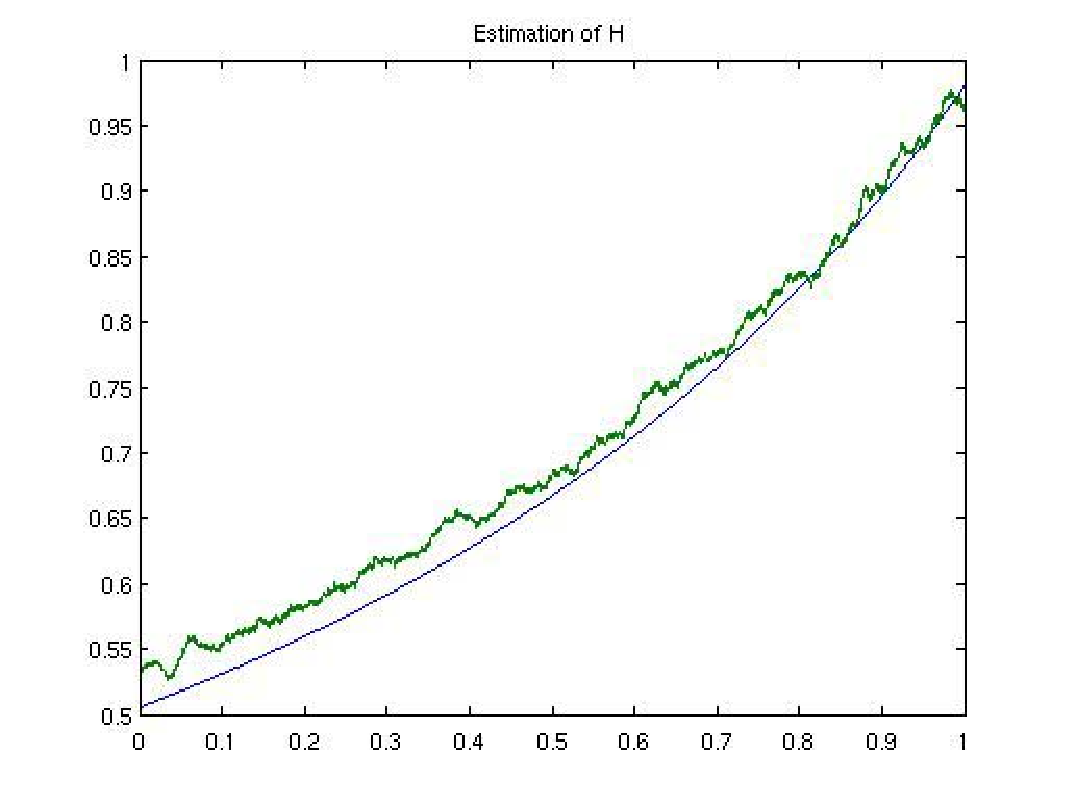} &
  \includegraphics[scale=0.3]{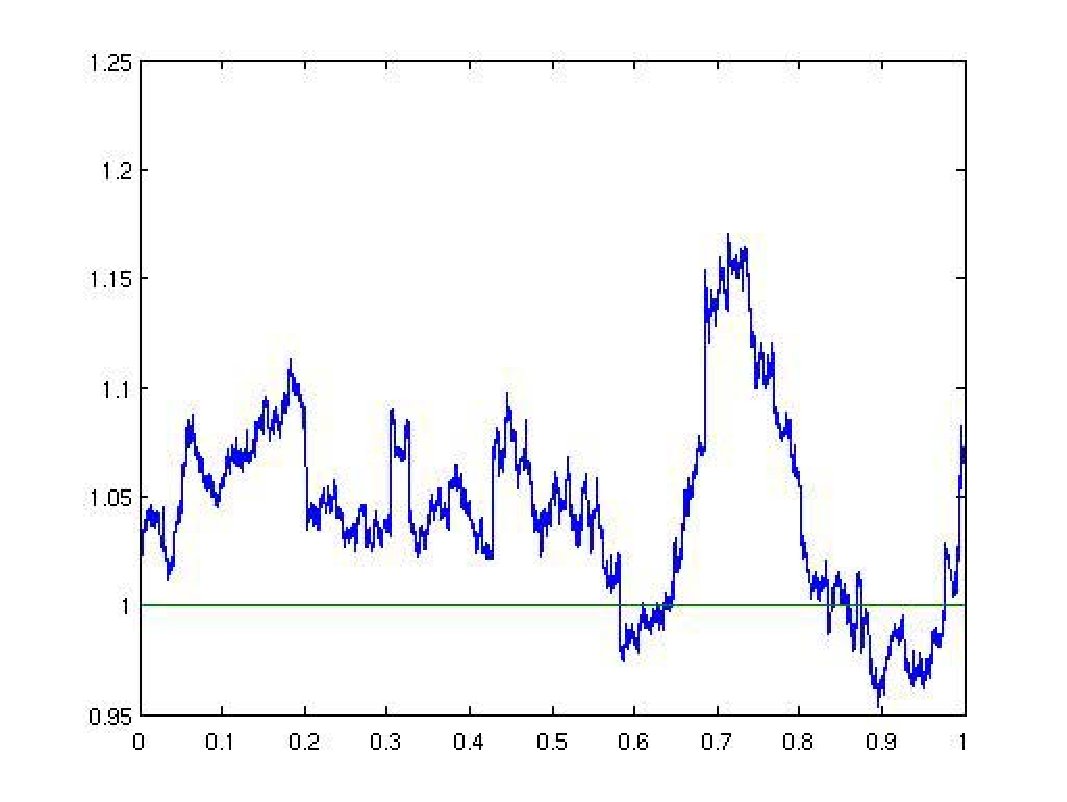}\\
  {\scriptsize  \hspace{1.2cm }$\alpha(t)=1.98-0.96 t$} & {\scriptsize  \hspace{1.2cm }$H(t)=\dfrac{1}{1.98-0.96 t}$} & \\
   \includegraphics[scale=0.3]{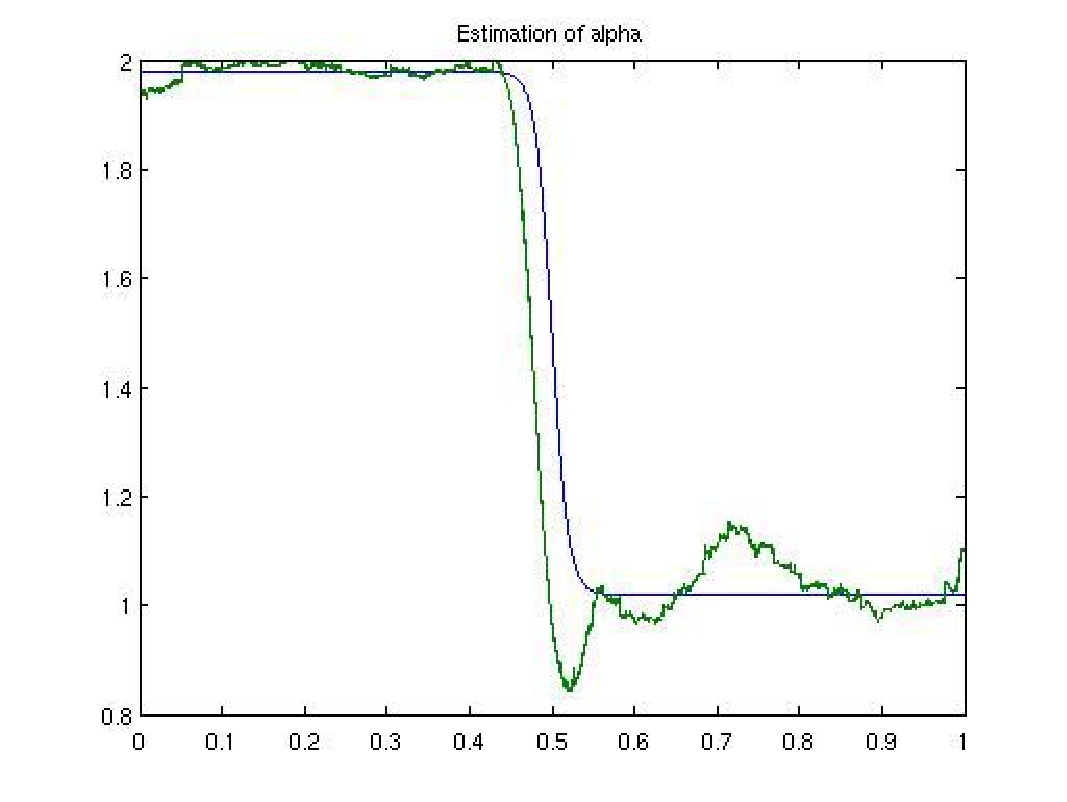} & \includegraphics[scale=0.3]{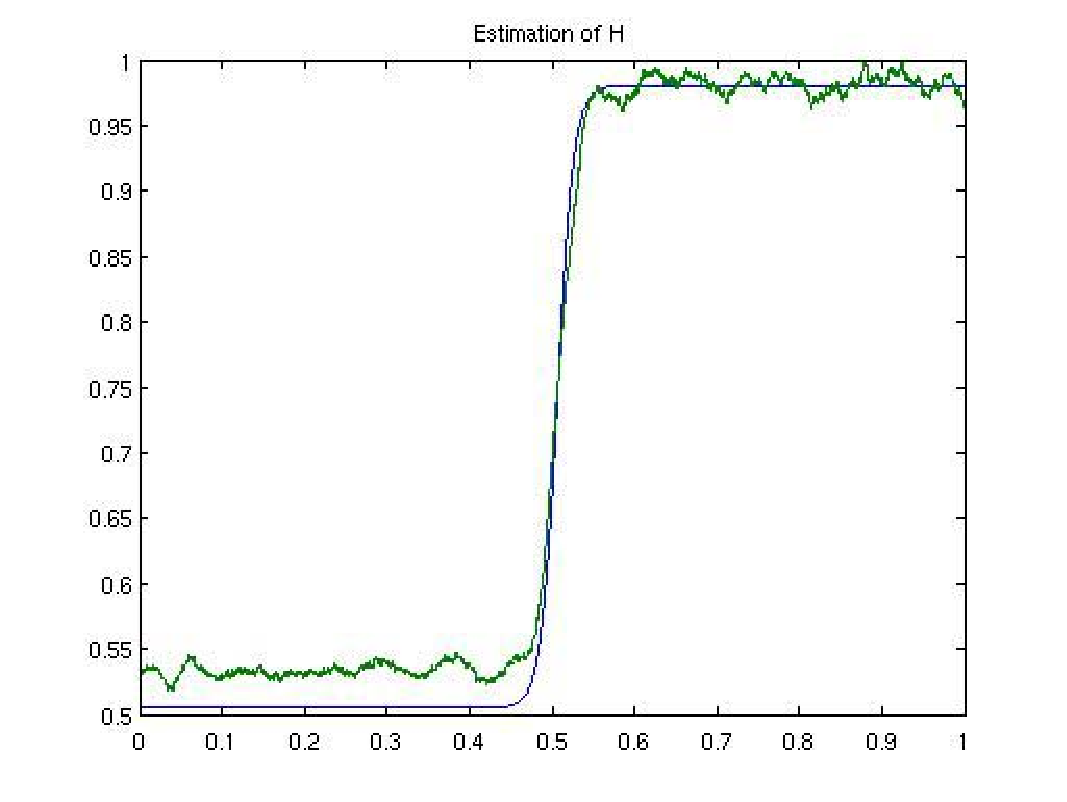} &
  \includegraphics[scale=0.3]{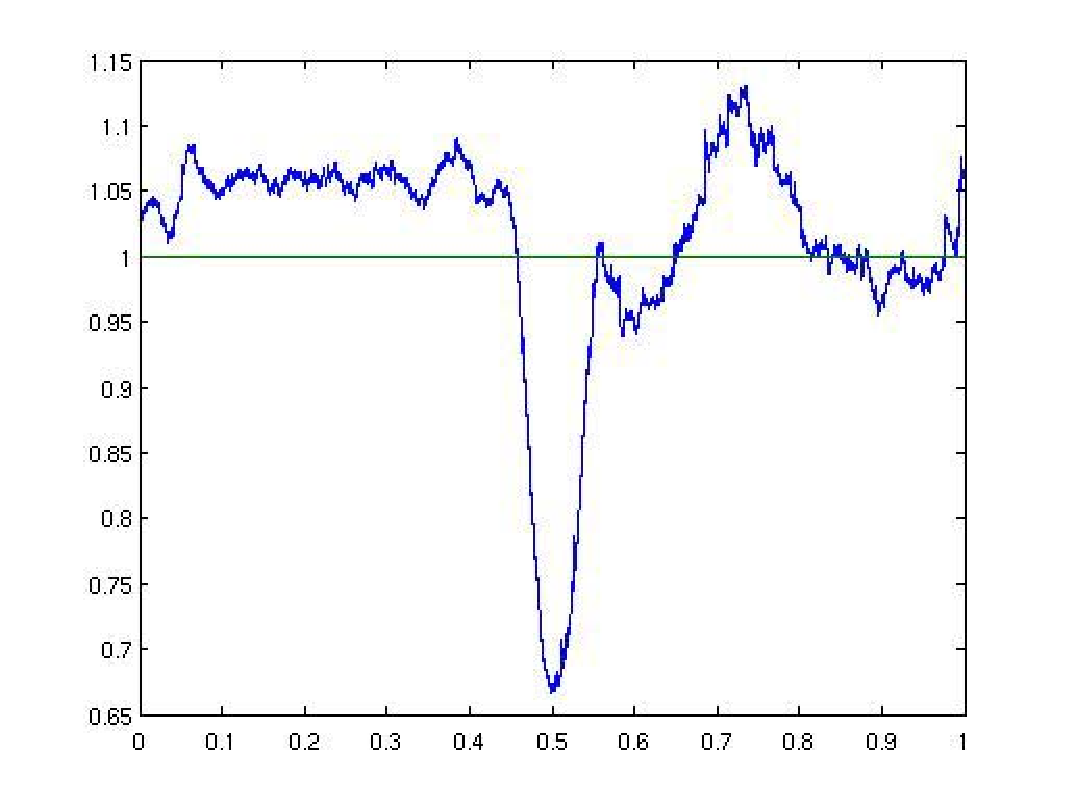}\\
  {\scriptsize  \hspace{0.7cm }$\alpha(t)=1.98-\dfrac{0.96}{1+\exp(20-40t)}$} & {\scriptsize  \hspace{0.7cm }$H(t)=\dfrac{1+\exp(20-40t)}{1.02+1.98\exp(20-40t)}$} & \\  
   \includegraphics[scale=0.3]{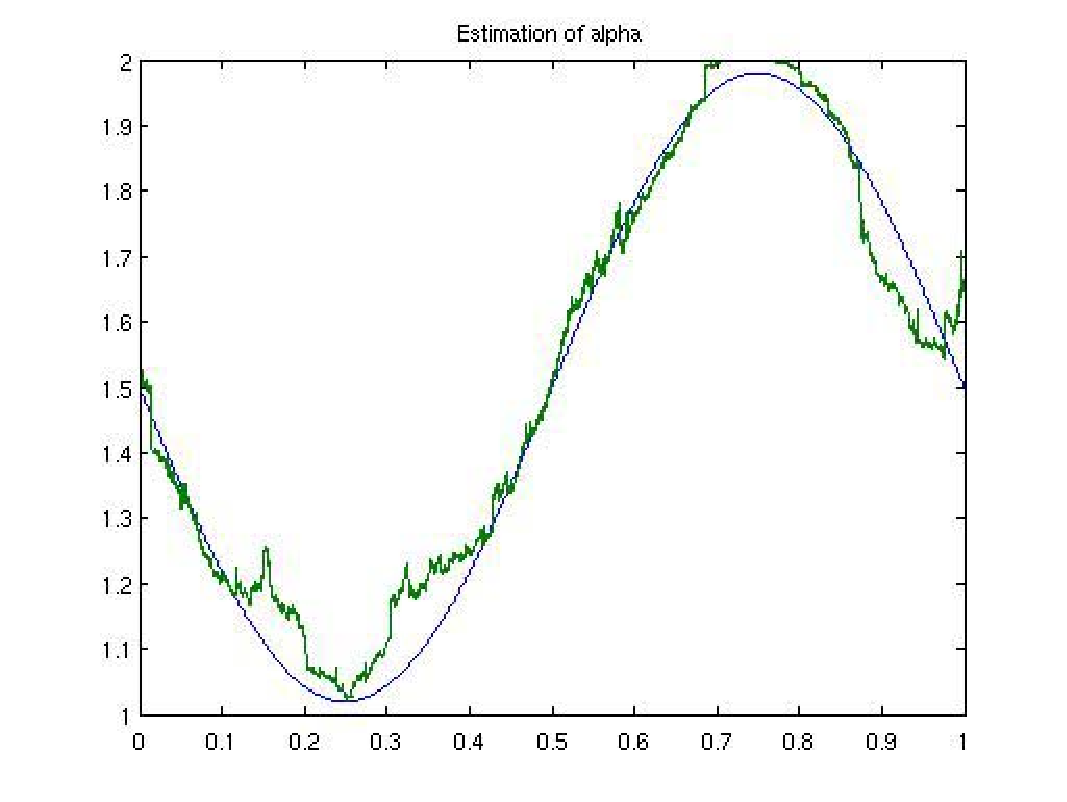} & \includegraphics[scale=0.3]{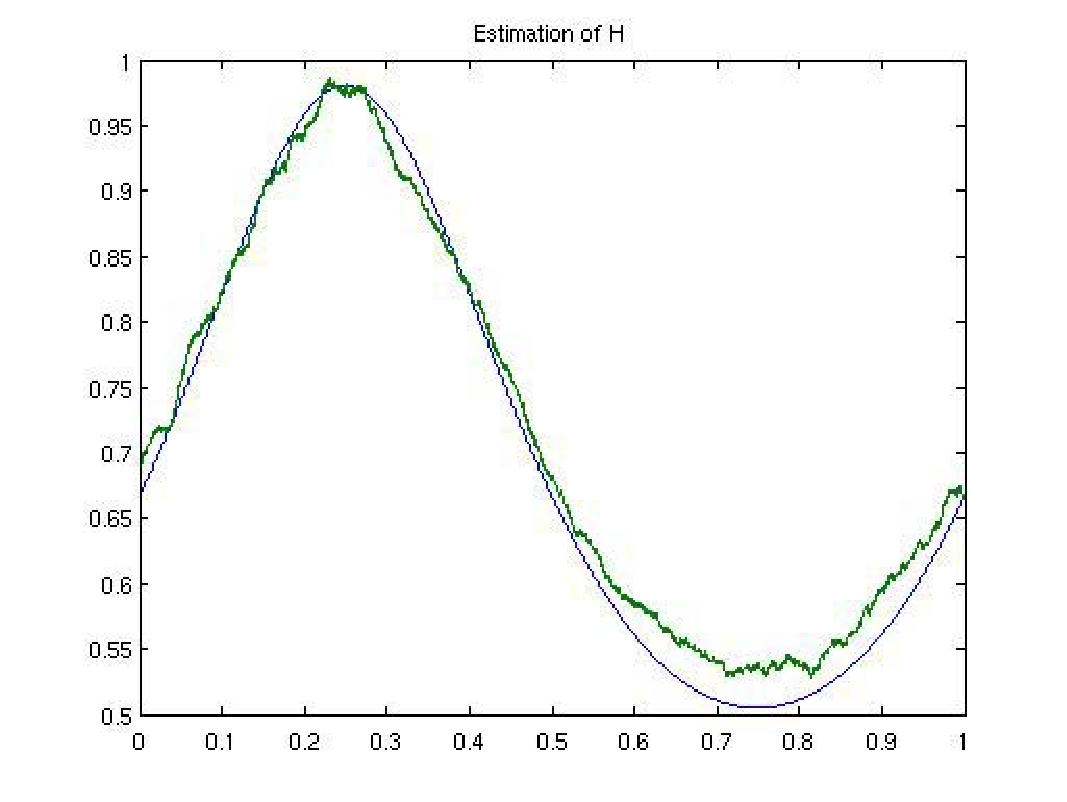} &
  \includegraphics[scale=0.3]{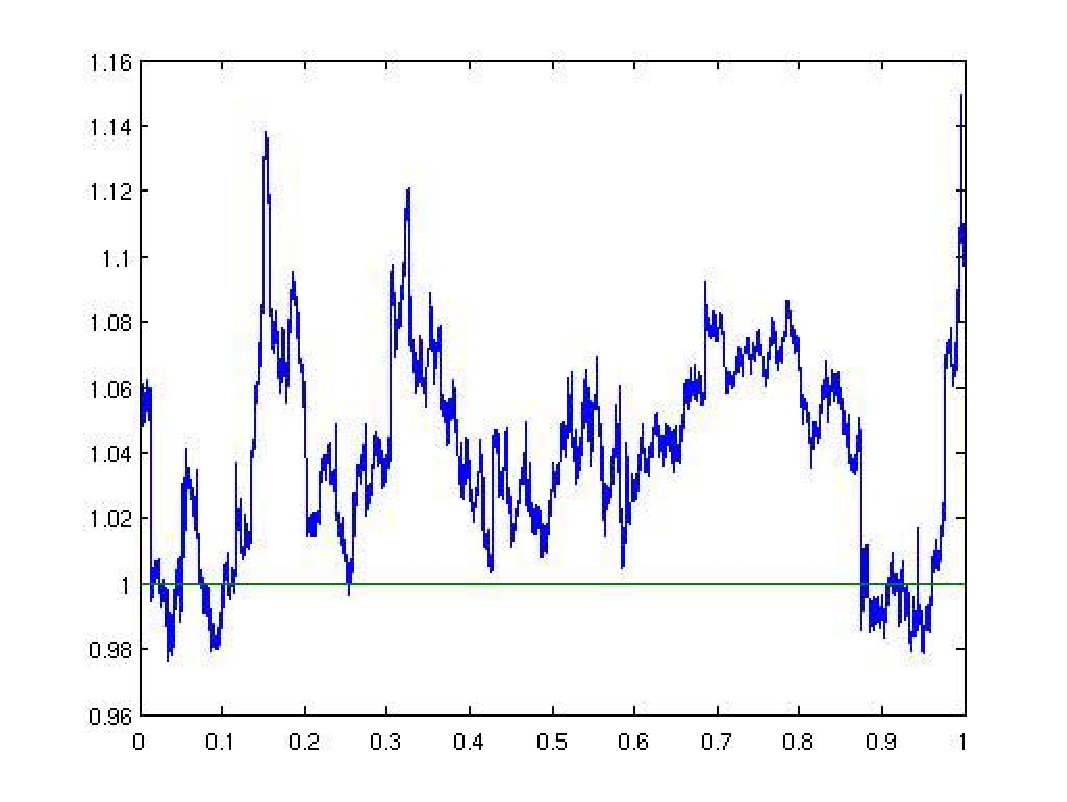}\\
  {\scriptsize  \hspace{0.9cm }$\alpha(t)=1.5-0.48\sin(2 \pi t)$} & {\scriptsize  \hspace{0.9cm }$H(t)=\dfrac{1}{1.5-0.48\sin(2 \pi t)}$} & \\
\end{tabular}
\caption{Trajectories on $(0,1)$ with $N=20000$ points, $n(N)=2042$ points for the estimator $\hat{\alpha}$, and $n(N)=500$ for $\hat{H}$. $\alpha$ and $\hat{\alpha}$ are represented in the first column, 
$H$ and $\hat{H}$ in the second column, and in the last column, we have drawn the product $ \hat{\alpha} \hat{H}$.}\label{fig1}
\end{figure}

We display on Figure \ref{fig1} some examples of estimations for various functions $\alpha$, the function $H$ satisfying the relation $H(t)=\frac{1}{\alpha(t)}$. The trajectories have been simulated using the field (\ref{fieldLevy}). For each $u \in (0,1)$, $X(.,u)$ is a $\alpha(u)$-stable L\'evy Motion. 
It is then an $\alpha(u)$-stable process with independent increments. We have generated these increments using the RSTAB program available in \cite{ST3} or in \cite{ST}, and then taken the diagonal $X(t,t)$.


Each function is pretty well-evaluated. We are able to recreate with the estimators the shape of the functions. However, we notice a significant bias on Figure \ref{fig1} in the estimation of $H$. It seems to decrease when $H$ is getting values close to 1. 
We observe this phenomenon with most trajectories, while the estimator $\hat{\alpha}$ seems to be unbiased. We have displayed the product $\hat{\alpha}\hat{H}$ in order to show the link between the estimators. We actually find again the asymtpotic relationship $H(t)=\frac{1}{\alpha(t)}$.


\begin{figure}[H]
\begin{tabular}{lll}
  \includegraphics[scale=0.3]{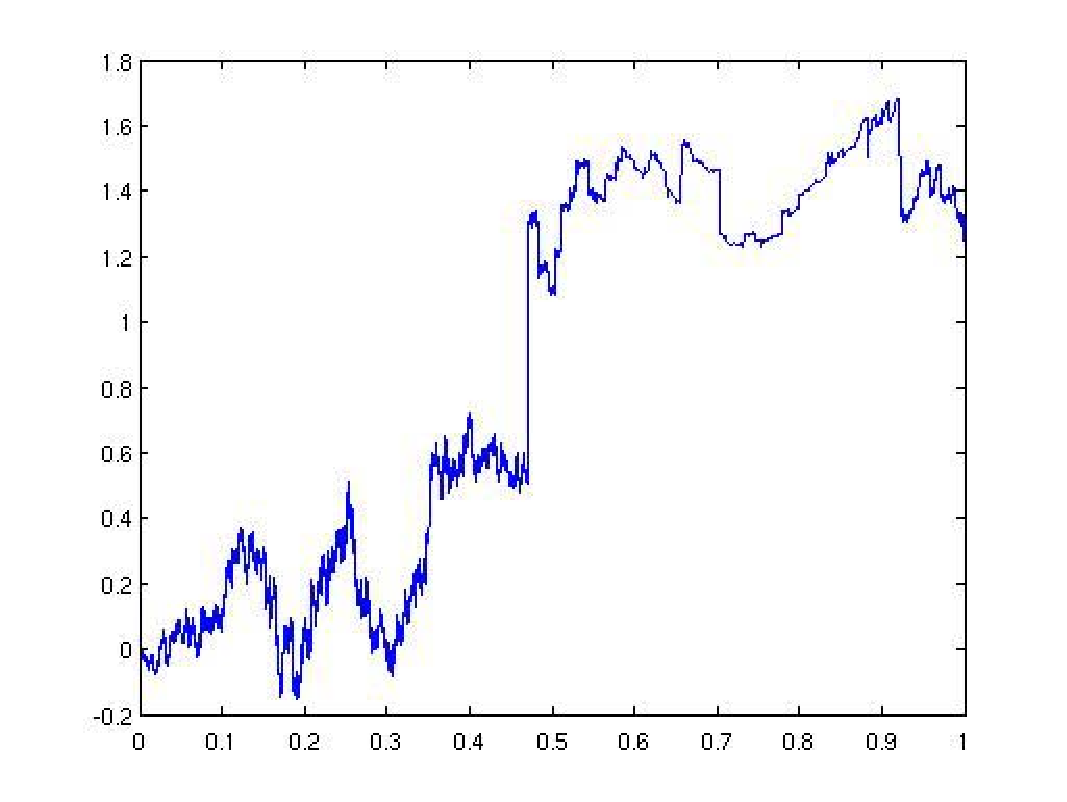} & \includegraphics[scale=0.3]{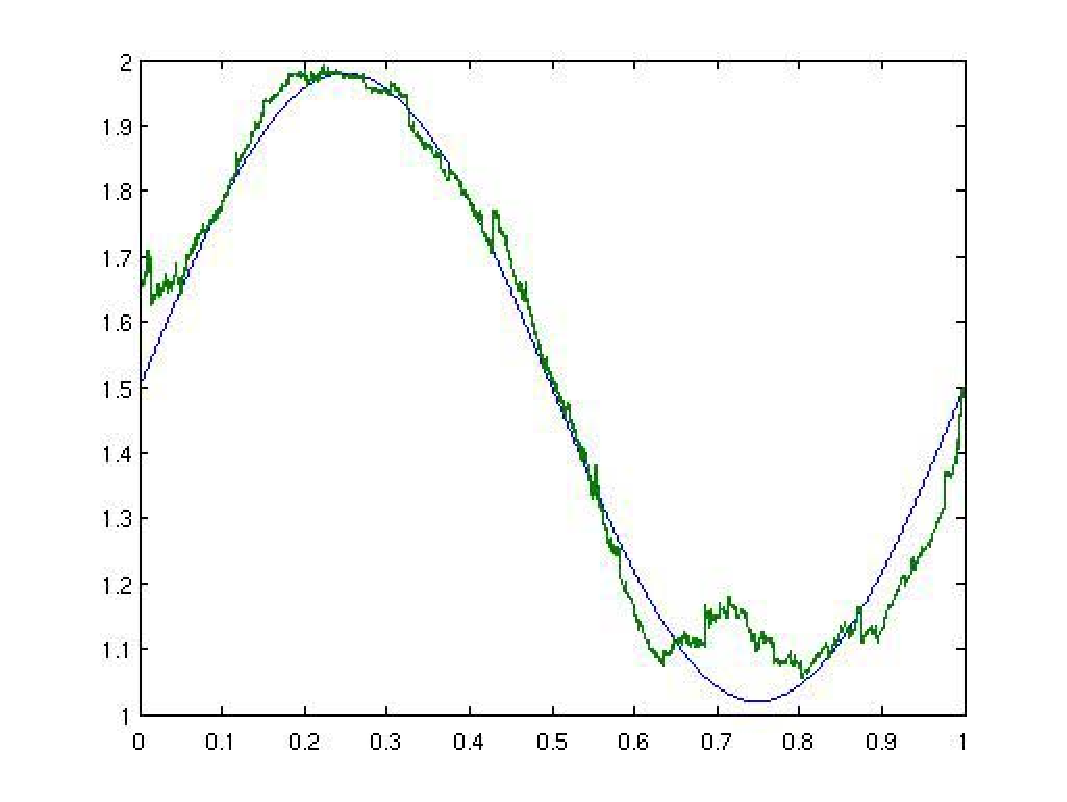} &
  \includegraphics[scale=0.3]{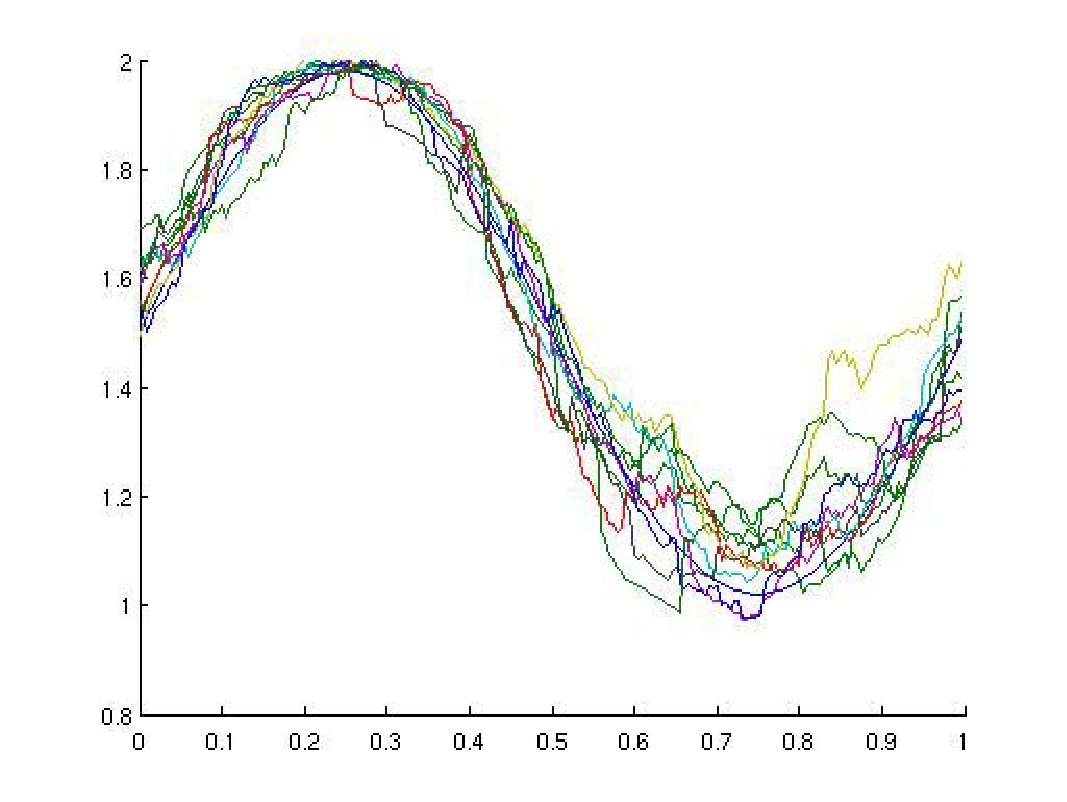}\\
  {\footnotesize \hspace{2.3cm}  a)} & {\footnotesize \hspace{2.3cm}  b)} & {\footnotesize \hspace{2.3cm}  c)}\\
\end{tabular}
\caption{Trajectory of a Levy process with $\alpha(t)= 1.5+0.48\sin(2 \pi t)$ in figure a), and the corresponding estimation of $\alpha$ in figure b) with $n(N)=2042$. The figure c) represents various estimations of $\alpha$ for the same function $\alpha(t)= 1.5+0.48\sin(2 \pi t)$, with different trajectories.}\label{fig2}
\end{figure}
We observe on Figure \ref{fig2} an evolution of the variance in the estimation of $\alpha$. It seems to increase when the function $\alpha$ is decreasing, and we conjecture that the variance at the point $t_0$ depends on the value $\alpha(t_0)$ in this way. 
In fact, the increments $Y_{k,N}$ are asymptotically distributed as an $\alpha(t_0)$-stable variable, so we expect that $S_N$ and $R_{\textrm{exp}}^{(N)}$ have a variance increasing when $\alpha$ is decreasing.

%
%
  
\begin{figure}[H]
\begin{tabular}{ll}
  \includegraphics[scale=0.3]{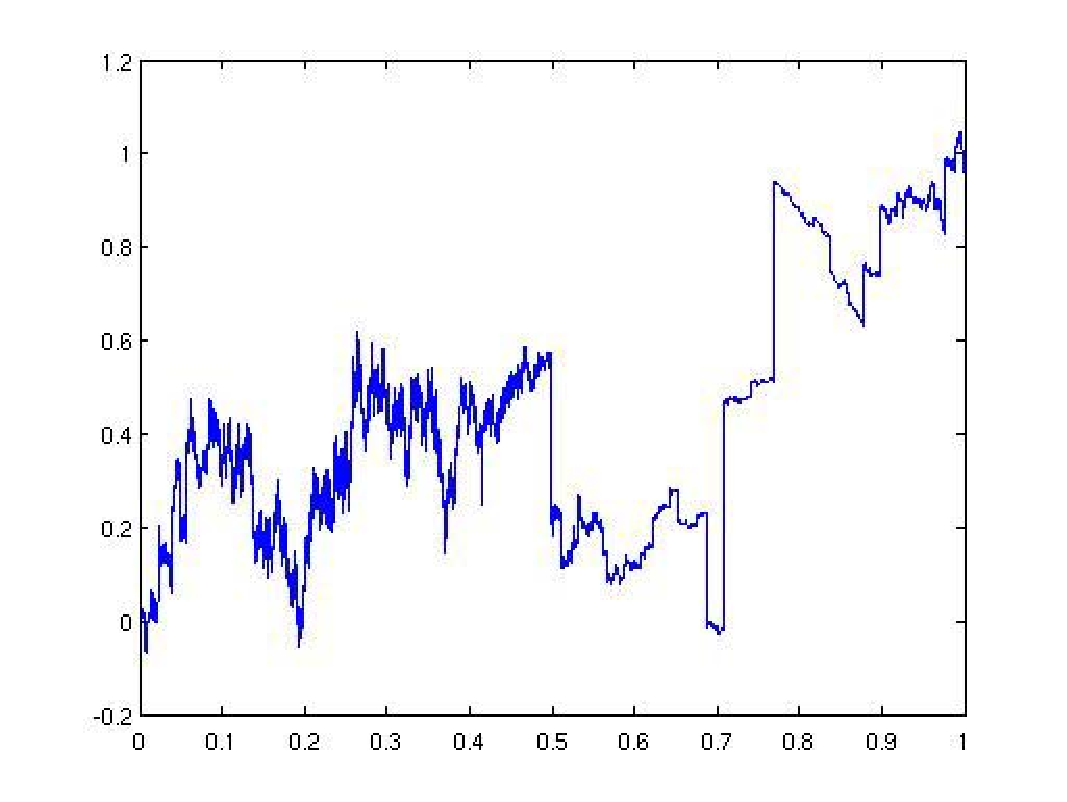} & \includegraphics[scale=0.3]{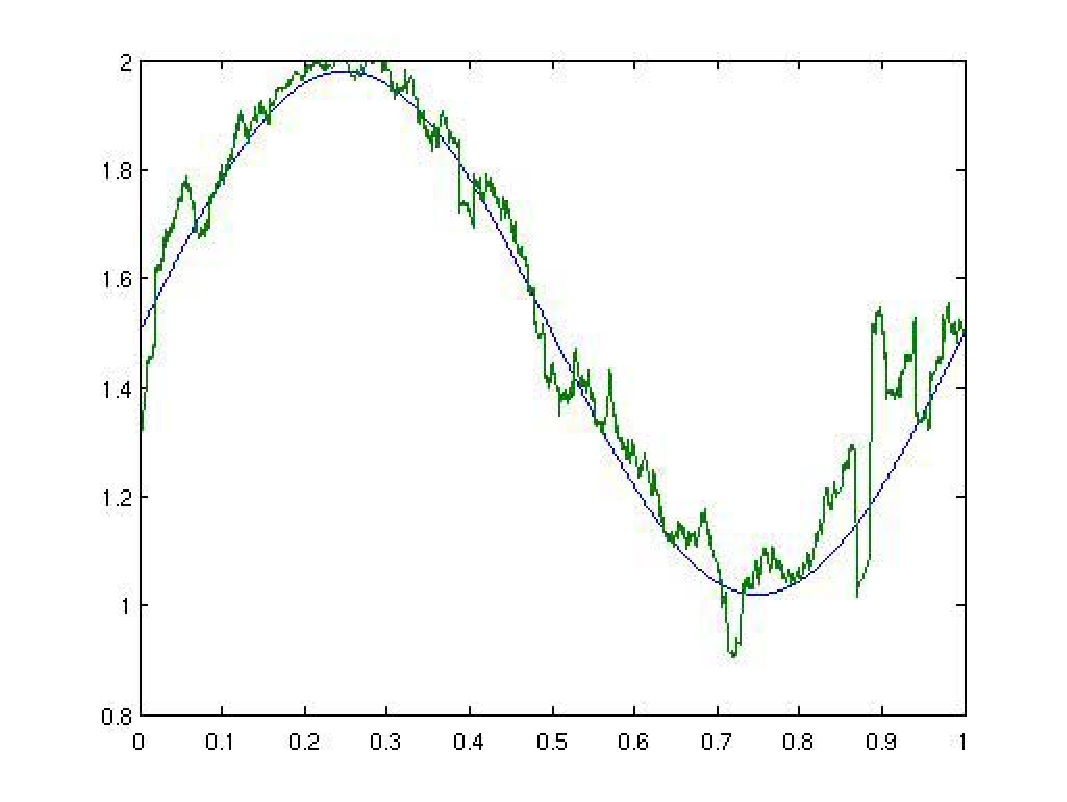} \\
  {\footnotesize \hspace{2.3cm}  d)} & {\footnotesize \hspace{2.3cm}  e)}\\
\end{tabular}  
\caption{Trajectory with $N=200000$ in figure d), and the estimation with $n(N)=3546$ in figure e).}\label{fig3}
\end{figure}

We have increased the resolution on Figure \ref{fig3}, taking more points for the discretization. The distance  observed on Figure \ref{fig2}.b for $\alpha$ near $1$ is then corrected.

 \subsection{Simulations with electrocardiogram}\label{AppliECG}
 
We consider an example of trajectory with a varying index of stability and a varying index of localisability.  The dataset comes from \cite{PH}.
 
\begin{figure}[H]
 \begin{center}
\begin{tabular}{l}
   \includegraphics[height=4cm, width=13cm]{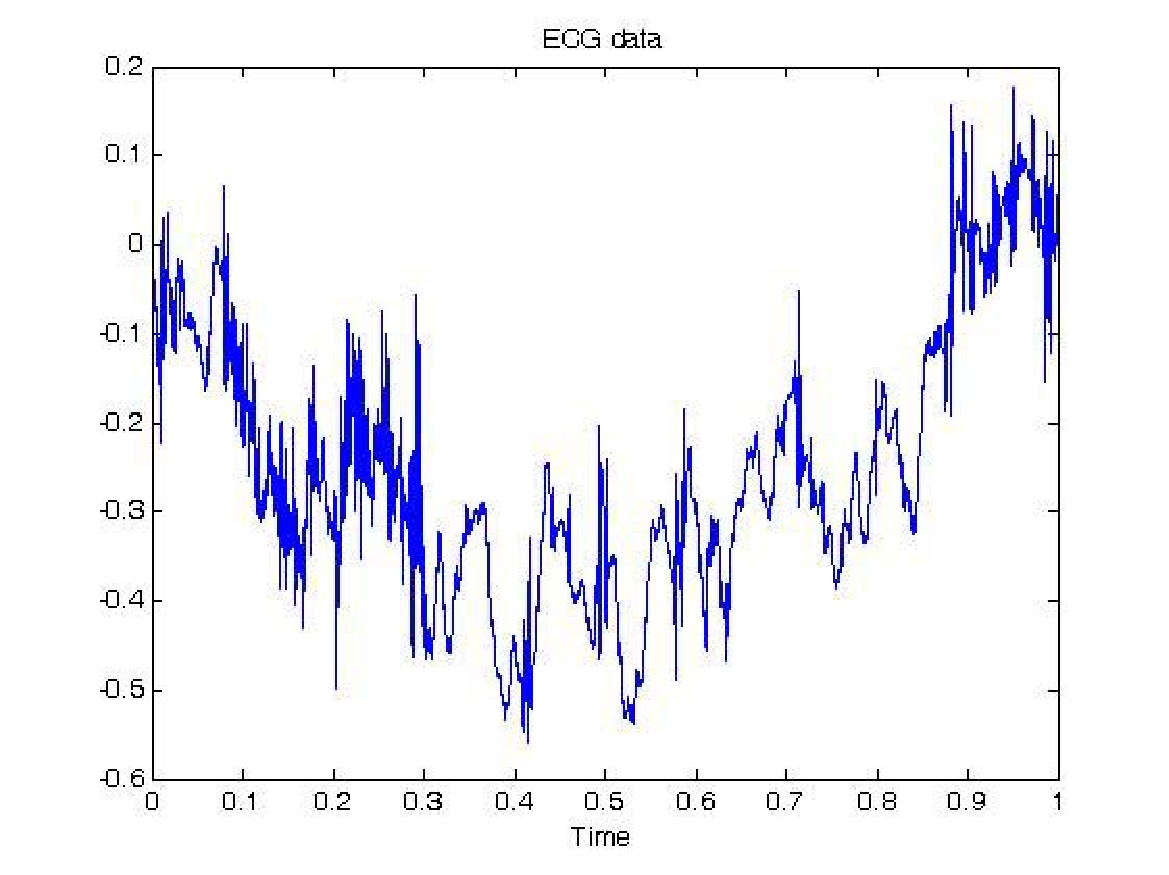} \\
\end{tabular}
\end{center}
\caption{Trajectory of the process $Y$ associated to ECG series with $N=1000000$ points.}\label{EcgTraj}
\end{figure}

 We denote $Z$ the process corresponding to an electrocardiogram. Its length is $N=1000000$ points. We consider then the process $Y$ defined by
 $$Y(j)=\sum_{i=1}^{j} \left(Z(i) - \frac{1}{N}\sum_{k=1}^{N}Z(k) \right).$$
 The realization of process $Y$ associated to EGC series is represented in Figure \ref{EcgTraj}. The increments of this process can not be regarded as stationary. We see in this example that the smoothness, as the intensity of significant jumps, is actually varying with time.
 
 
\begin{figure}[H]
 \begin{center}
\begin{tabular}{l}
   \includegraphics[height=4cm, width=13cm]{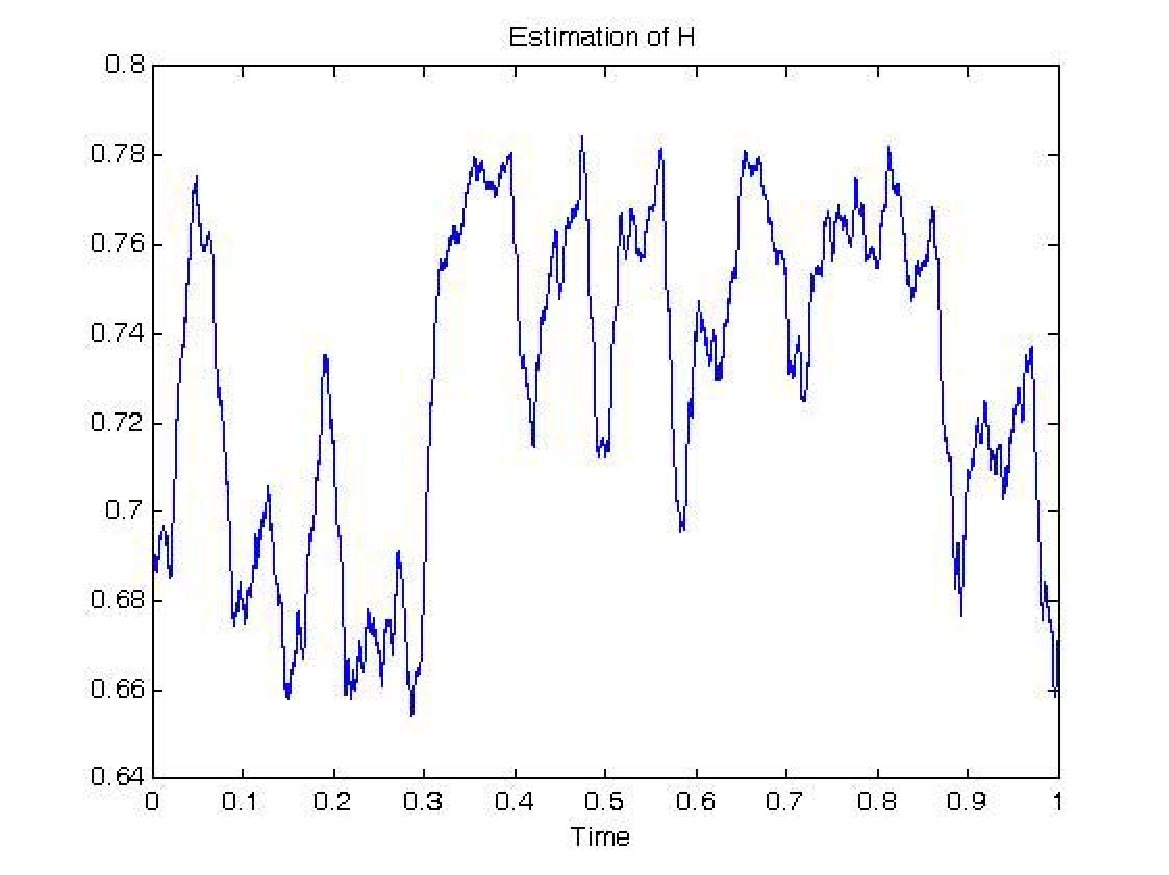} \\
\end{tabular}
\end{center}
\caption{Estimation of $H$ calculated for the process represented in Figure \ref{EcgTraj}.}\label{HEcgTraj}
\end{figure}

  We have done an estimation of the localisability function $H$ for this process $Y$. Figure \ref{HEcgTraj} represents an estimation of $H$ as function of $t$. The estimate of $H$ is calculated by taking $n(N)=25000$ points.

\begin{figure}[H]
 \begin{center}
\begin{tabular}{l}
   \includegraphics[height=4cm, width=13cm]{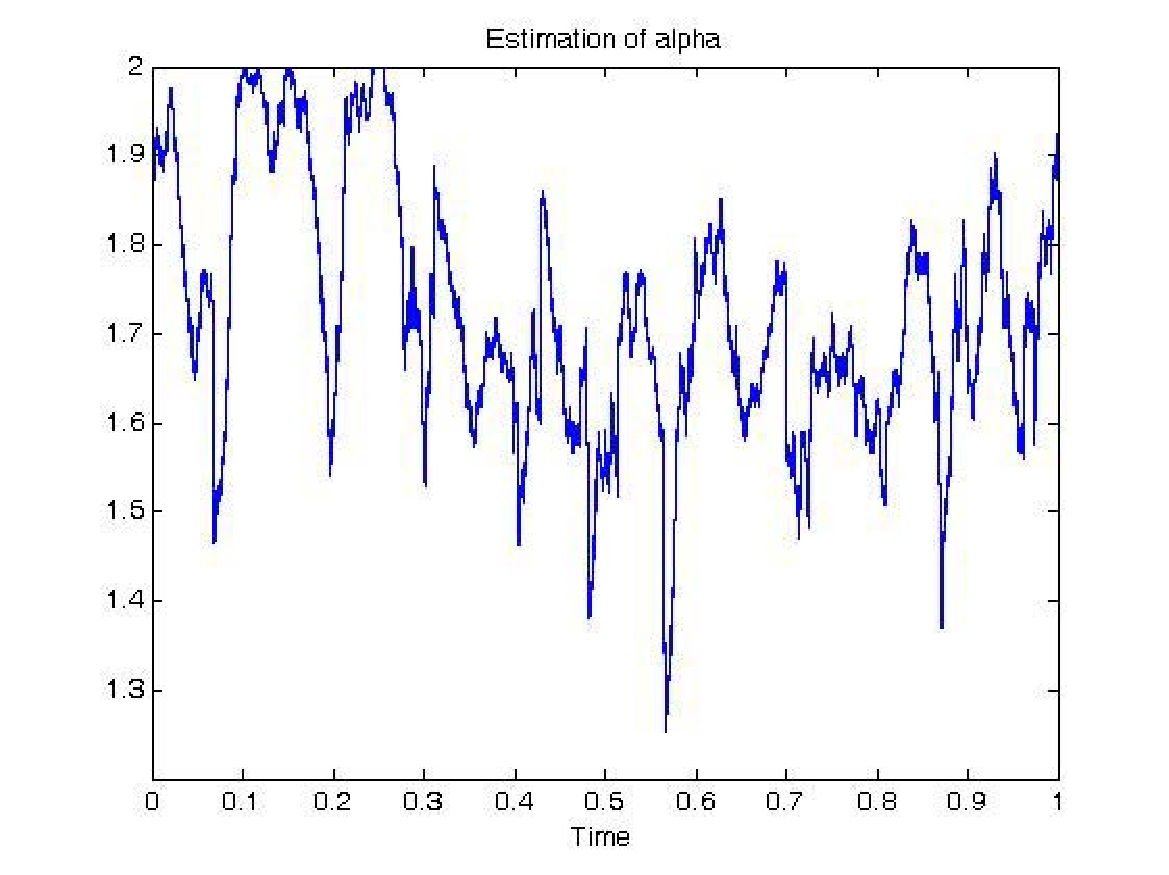} \\
\end{tabular}
\end{center}
\caption{Estimation of $\alpha$ calculated for the process of Figure \ref{EcgTraj}.}\label{AlphaEcgTraj}
\end{figure}

We notice a correlation between the noisy areas of the trajectory and the times when the exponent $H$ is small, and also a greatest exponent when the trajectory seems to be smoother. For the estimation of the function $\alpha$, we have taken $n(N)=25000$ too. The result is presented in Figure \ref{AlphaEcgTraj}.
 We observe also here a link between the noise and the function $\alpha$. When the intensity of the significant jumps of the trajectory is high, the stability function is close to $2$. A lower stability index matches to a period with a lower intensity of significant jumps.


\section{Assumptions}\label{Assum}

This section gathers the various conditions required on the considered processes so that
our results hold. These asumptions are of three kinds: regularity condition that entail localisability, moment conditions related to the fact that we work in certain functional spaces and finally, H\"{o}lder conditions which enable to transfer the behaviour of $f$ to the one of $Y$.

{\bf Regularity}
 
\begin{itemize}
    \item (R1) The family of functions $v \to f(t,v,x)$ is differentiable for all $(v,t)$ in $U^2$ and almost all $x$ in $E$. The derivatives of $f$ with respect to $v$ are denoted by $f'_v$.
\end{itemize}

{\bf Moments conditions}

\begin{itemize}    
    
    \item (M1) There exists $\delta > \frac{d}{c} - 1$ such that :
\begin{displaymath}\label{kercond2sf}
\sup_{t \in U}  \int_\bbbr \left[ \sup_{w \in U} (|f(t,w,x)|^{\alpha(w)}) \right]^{1+\delta}  r(x)^{\delta} \hspace{0.1cm}  m(dx) < \infty.
\end{displaymath}
	
    \item (M2) There exists $\delta > \frac{d}{c} - 1$ such that :
\begin{displaymath}\label{kercond3sf}
\sup_{t \in U}  \int_\bbbr \left[ \sup_{w \in U} (|f'_v(t,w,x)|^{\alpha(w)}) \right]^{1+\delta} r(x)^{\delta} \hspace{0.1cm} m(dx) < \infty.
\end{displaymath}

    \item (M3) There exists $\delta > \frac{d}{c} - 1$ such that :
\begin{displaymath}\label{kercond5sf}
\sup_{t \in U}  \int_\bbbr \left[ \sup_{w \in U} \left[  \left|f(t,w,x) \log(r(x))  \right|^{\alpha(w)}  \right]\right]^{1+\delta} r(x)^{\delta} \hspace{0.1cm} m(dx) < \infty.
\end{displaymath}

   \item (M4) There exists $K_U >0$ such that $\forall v \in U$, $\forall u \in U$, $\forall x \in \bbbr$,
    \begin{displaymath}
    \left| f(v,u,x) \right| \leq K_U.
    \end{displaymath}

  \item (M5) There exists $K_U >0$ such that $\forall v \in U$, $\forall u \in U$, $\forall x \in \bbbr$,
    \begin{displaymath}
    \left| f'_v(v,u,x) \right| \leq K_U.
    \end{displaymath}

	\item (M6) There exists $ K_U >0$ such that $\forall  v \in U$, $\forall  u \in U$, 
	\begin{displaymath}
	\int_{\bbbr} \left| f(v,u,x) \right|^2 m(dx) \leq K_U.
	\end{displaymath}

    \item (M7) 
    \begin{displaymath}
    \inf_{v \in U} \int_{\bbbr} f(v,v,x)^2 m(dx) >0.
    \end{displaymath}
\end{itemize}

{\bf H\"{o}lder conditions}

\begin{itemize}

     \item (H1) There exists $K_U >0$ such that $\forall (u,v) \in U^2$, $\forall x \in \bbbr$,
    \begin{displaymath}
    \frac{1}{|v-u|^{H(u)-1/ \alpha(u)}}\left| f(v,u,x) -f(u,u,x) \right| \leq K_U.
    \end{displaymath}

	\item (H2) There exists  $K_U >0$ such that  $\forall (u,v) \in U^2$,
\begin{displaymath}	
	\frac{1}{|v-u|^{H(u)\alpha(u)}} \int_{\bbbr} | f(v,u,x)-f(u,u,x)|^{\alpha(u)} m(dx)  \leq K_U.
\end{displaymath}

	\item (H3)  There exists $ p \in (d,2)$, $p \geq 1$ and $ K_U >0$ such that $\forall (u,v) \in U^2$,
	\begin{displaymath}
	 \frac{1}{|v-u|^{1+p(H(u)-\frac{1}{\alpha(u)})}} \int_{\bbbr} \left| f(v,u,x) - f(u,u,x) \right|^p m(dx)\leq K_U.
	\end{displaymath}

	\item (H4) There exists a positive function $g$ defined on $U$ such that 
    \begin{displaymath}
\lim_{r \rightarrow 0} \sup_{t \in U} \left| \frac{1}{r^{1+2(H(t)-1/\alpha(t))})} \int_{\bbbr} \left( f(t+r,t,x) - f(t,t,x) \right)^2  m(dx) - g(t) \right|  = 0.
\end{displaymath}

	\item (H5) There exists $K_U >0$ such that $\forall (u,v) \in U^2$, 
	\begin{displaymath}
	 \frac{1}{|v-u|^2} \int_{\bbbr} \left| f(v,v,x) - f(v,u,x) \right|^2 m(dx)\leq K_U.
	\end{displaymath}

\end{itemize}

\section{Proofs}\label{proo}

In all the proofs, $K_U$ denotes a generic constant which depends on the interval $U$ and may vary from line to line.

{\bf Proof of Lemma \ref{LemSup}}

Let $B \in \bbbr$, $B \geq \max(5, \frac{6}{c} )$. Let $[a,b] \subset U$. We denote $E_N = \{ k \in \bbbn \cap [0,N-1], \frac{k}{N} \in [a,b] \textrm{ or } \frac{k+1}{N} \in [a,b] \}.$
For $N$ large enough, since $\lim_{N \rightarrow +\infty}\limits \frac{n(N)}{N} = 0$, for all $k \in E_N$ and $j \in \bbbn$ such that $k-\frac{n(N)}{2} \leq j \leq  k+\frac{n(N)}{2}-1$, $\frac{j}{N} \in U$ and $\frac{j+1}{N} \in U$.
The function $t \mapsto \hat{H}_N(t)$ is a step function so

\begin{eqnarray*}
 \P \left( \sup_{t \in [a,b]}\limits | \hat{H}_N(t)| > B \right) & \leq & \P \left( \cup_{k \in E_N} \{ |\sum_{j=k-\frac{n(N)}{2}}^{k+\frac{n(N)}{2}-1}\limits \log | Y_{j,N}| | > B n(N) \log N \} \right)\\
 & \leq & \sum_{k \in E_N} \sum_{j=k-\frac{n(N)}{2}}^{k+\frac{n(N)}{2}-1}\limits  \P \left( | \log | Y_{j,N}| | > B \log N \right)\\
 & \leq & \sum_{k \in E_N} \sum_{j=k-\frac{n(N)}{2}}^{k+\frac{n(N)}{2}-1}\limits  \P \left( |X(\frac{j+1}{N},\frac{j+1}{N}) - X(\frac{j+1}{N},\frac{j}{N}) | \geq \frac{N^B}{2} \right)\\
 & & + \sum_{k \in E_N} \sum_{j=k-\frac{n(N)}{2}}^{k+\frac{n(N)}{2}-1}\limits \P \left( |X(\frac{j+1}{N},\frac{j}{N}) - X(\frac{j}{N},\frac{j}{N}) | \geq \frac{N^B}{2}  \right) \\
 & & + \sum_{k \in E_N} \sum_{j=k-\frac{n(N)}{2}}^{k+\frac{n(N)}{2}-1}\limits \P \left( |Y(\frac{j+1}{N}) - Y(\frac{j}{N}) | \leq \frac{1}{N^B}\right). \\
 \end{eqnarray*}
We control each probability of the right term. With the conditions {\bf (R1)}, {\bf (M1)}, {\bf (M2)} and {\bf (M3)}, we can apply Proposition 4.9 of \cite{LGLV3} : there exists $K_U > 0$ such that for all $(u,v) \in U^2$ and $x > 0$, 
\begin{equation}\label{majprob}
\P \left(|X(v,v)-X(v,u)| > x \right) \leq K_U \left( \frac{|v-u|^d}{x^d}(1 +|\log \frac{|v-u |}{x}|^d) + \frac{|v-u|^c}{x^c}(1+|\log \frac{|v-u |}{x}|^c) \right).
\end{equation}

We obtain the existence of a constant $K>0$ which depends on $U$, $B$, $c$ and $d$ such that 

\begin{displaymath}
 \P \left( |X(\frac{j+1}{N},\frac{j+1}{N}) - X(\frac{j+1}{N},\frac{j}{N}) | \geq \frac{N^B}{2} \right) \leq K \frac{|\log N|^d}{N^{Bc}}.
\end{displaymath}

The process $X(.,\frac{j}{N})$ is an $\alpha(\frac{j}{N})$-stable process, so 
\begin{eqnarray*}
 \P \left( |X(\frac{j+1}{N},\frac{j}{N}) - X(\frac{j}{N},\frac{j}{N}) | \geq \frac{N^B}{2}  \right) &\leq &\frac{2^{c/2}}{N^{\frac{Bc}{2}}} \E \left[ | X(\frac{j+1}{N},\frac{j}{N}) - X(\frac{j}{N},\frac{j}{N}) |^{c/2}\right]\\
 & = & \frac{K_1}{N^{\frac{Bc}{2}}} \left[ \int_E |f(\frac{j+1}{N},\frac{j}{N},x) - f(\frac{j}{N},\frac{j}{N},x) |^{\alpha(\frac{j}{N})} m(dx) \right]^{\frac{c}{2 \alpha(\frac{j}{N})}} \\
\end{eqnarray*}
where $K_1 = \frac{2^{c} \Gamma(1-\frac{c}{2\alpha(\frac{j}{N})})}{c \int_{0}^{+\infty} u^{-\frac{c}{2}-1} \sin^2(u) du}$. With the condition {\bf (H2)}, we obtain a constant $K_U >0$ such that $  \P \left( |X(\frac{j+1}{N},\frac{j}{N}) - X(\frac{j}{N},\frac{j}{N}) | \geq \frac{N^B}{2}  \right) \leq \frac{K_U}{N^{\frac{Bc}{2}}}$. 
With the conditions {\bf (R1)}, {\bf (M4)}, {\bf (M5)}, {\bf (M6)}, {\bf (M7)}, {\bf (H1)}, {\bf (H3)}, {\bf (H4)} and {\bf (H5)}, we use for the third term Propositions 4.10 and 4.8 of \cite{LGLV3}: there exists $K>0$ and $N_0 \in \bbbn$ such that for all $t \in U$, for all $N \geq N_0$ and all $x >0$,
\begin{equation}\label{probinferi}
 \P \left( | Y(t+\frac{1}{N}) - Y(t) | < x\right) \leq K N^{H(t)} x.
\end{equation}
Then 
\begin{equation}
 \P \left( |Y(\frac{j+1}{N}) - Y(\frac{j}{N}) | \leq \frac{1}{N^B}\right) \leq \frac{K}{N^B} N^{H(\frac{j}{N})} .
\end{equation}
We get then 
\begin{displaymath}
 \P \left( \sup_{t \in [a,b]}\limits | \hat{H}_N(t)| > B \right) \leq K_U N n(N) \left( \frac{|\log N|^d}{N^{Bc}} + \frac{1}{N^{\frac{Bc}{2}}} + \frac{1}{N^{B-H_+}} \right),
\end{displaymath}
and we conclude with the Borel Cantelli lemma \Box

{\bf Proof of Theorem \ref{ConvLpSnp}}

First, note that the condition $\lim_{j \rightarrow +\infty}\limits \int_E | h_{0,t_0}(x) h_{j,t_0}(x)|^{\frac{\alpha(t_0)}{2}} m(dx)=0$ implies the following condition:
\begin{itemize}
 \item (C*)  There exists $\varepsilon_1 >0$ and $j_0 \in \mathbb{N}$ such that for all $j\geq j_0$,
 \begin{displaymath}
  \int_E | h_{0,t_0}(x) h_{j,t_0}(x)|^{\frac{\alpha(t_0)}{2}} m(dx) \leq (1-\varepsilon_1) \Vert h_{0,t_0} \Vert_{\alpha(t_0)}^{\alpha(t_0)},
 \end{displaymath}
\end{itemize}

Let $p \in [p_0,\alpha(t_0) )$. We define
\begin{displaymath}
 A_N(p) = \frac{N^{pH(t_0)}}{n(N)} \sum_{k=[Nt_0]-\frac{n(N)}{2}}^{[Nt_0]+\frac{n(N)}{2}-1} \left| X(\frac{k+1}{N},\frac{k+1}{N}) - X(\frac{k+1}{N},t_0)\right|^p,
\end{displaymath}
\begin{displaymath}
  B_N(p)=\frac{N^{pH(t_0)}}{n(N)} \sum_{k=[Nt_0]-\frac{n(N)}{2}}^{[Nt_0]+\frac{n(N)}{2}-1} \left| X(\frac{k}{N},\frac{k}{N}) - X(\frac{k}{N},t_0)\right|^p 
\end{displaymath}
and
\begin{displaymath}
 C_N(p) = \frac{N^{pH(t_0)}}{n(N)} \sum_{k=[Nt_0]-\frac{n(N)}{2}}^{[Nt_0]+\frac{n(N)}{2}-1} \left| X(\frac{k+1}{N},t_0) - X(\frac{k}{N},t_0)\right|^p.
\end{displaymath}

Let $Z=X(1,t_0).$ We have, for $p \leq 1$,
\begin{eqnarray*}
 \P \left( |N^{pH(t_0)}S_N^p(p) - \E|Z|^{p} | >x \right) & \leq & \P \left( |N^{pH(t_0)}S_N^p(p) - C_N(p) | \geq \frac{x}{2} \right) + \P \left( |\E|Z|^{p} - C_N(p) | \geq \frac{x}{2}\right)\\
& \leq & \P \left( |\E|Z|^{p} - C_N(p) | \geq \frac{x}{2}\right) + \P \left( A_N(p) + B_N(p) \geq \frac{x}{2}\right)\\
\end{eqnarray*}
and for $p\geq 1$,
\begin{eqnarray*}
 \P \left( |N^{H(t_0)}S_N(p) - (\E|Z|^{p})^{\frac{1}{p}} | >x \right) & \leq & \P \left( |N^{H(t_0)}S_N(p) - C_N^{\frac{1}{p}}(p) | \geq \frac{x}{2} \right) \\
 & &+ \P \left( |C_N^{\frac{1}{p}}(p) - (\E|Z|^{p})^{\frac{1}{p}}| \geq \frac{x}{2}\right)\\
 & \leq &  \P \left( |(\E|Z|^{p})^{\frac{1}{p}} - C_N^{\frac{1}{p}}(p) | \geq \frac{x}{2}\right) +\P \left( A_N^{\frac{1}{p}}(p) + B_N^{\frac{1}{p}}(p)\geq \frac{x}{2}\right) . \\
 \end{eqnarray*}

To prove Theorem \ref{ConvLpSnp}, it is enough to show that $A_N(p) \overset{\P}{\longrightarrow} 0$, $B_N(p) \overset{\P}{\longrightarrow} 0$ and $C_N(p) \overset{\P}{\longrightarrow} \E|Z|^{p} $.

We consider first $A_N(p) \overset{\P}{\longrightarrow} 0$. Let $\delta_N(dt) = \frac{N}{n(N)} \mathbf{1}_{\{ \frac{[Nt_0]}{N}-\frac{n(N)}{2N} \leq t < \frac{[Nt_0]}{N}+\frac{n(N)}{2N}\}} dt$. Let $U$ be an open interval satisfying the conditions of the theorem and $t_0 \in U$. 
We can fix $N_0 \in \mathbb{N}$ and $V \subset U$ an open interval depending on $t_0$ such that for all $N \geq N_0$ and all $t \in V$, $ \frac{[Nt]+1}{N} \in U$, $\frac{[Nt]}{N} \in U$, $\int_{0}^{1}  \delta_N(dt) =  \int_V \delta_N(dt)$, and such that the inequality (\ref{majprob}) holds.
\begin{eqnarray*}
 \P \left( A_N(p) > x \right) & = & \P \left( \int_{0}^{1} \left| \frac{X(\frac{[Nt]+1}{N},\frac{[Nt]+1}{N}) - X(\frac{[Nt]+1}{N},t_0)}{(1/N)^{H(t_0)}}\right|^p \delta_N(dt) > x \right) \\
 & \leq & \frac{1}{x} \int_V \E \left[\left| \frac{X(\frac{[Nt]+1}{N},\frac{[Nt]+1}{N}) - X(\frac{[Nt]+1}{N},t_0)}{(1/N)^{H(t_0)}}\right|^p\right] \delta_N(dt)\\
\end{eqnarray*}
Let $t \in V$.
\begin{displaymath}
\E \left[\left| \frac{X(\frac{[Nt]+1}{N},\frac{[Nt]+1}{N}) - X(\frac{[Nt]+1}{N},t_0)}{(1/N)^{H(t_0)}}\right|^p\right] = \int_{0}^{\infty} \P \left( \left| \frac{X(\frac{[Nt]+1}{N},\frac{[Nt]+1}{N}) - X(\frac{[Nt]+1}{N},t_0)}{(1/N)^{H(t_0)}} \right| > u^{1/p}\right) du.
\end{displaymath}
Let $u>0$. We know from (\ref{majprob}) that there exists $K_U >0$ such that for all $t \in V$,

\begin{displaymath}
\P \left( \left| \frac{X(\frac{[Nt]+1}{N},\frac{[Nt]+1}{N}) - X(\frac{[Nt]+1}{N},t_0)}{(1/N)^{H(t_0)}} \right| > u^{1/p}\right) \leq  K_U \frac{((\log N)^c+|\log u|^c )}{N^{c(1-H(t_0))}u^{c/p}}  + K_U \frac{((\log N)^d+|\log u|^d )}{N^{d(1-H(t_0))}u^{d/p}},
\end{displaymath}
so, with the assumption $H(t_0) <1$,
\begin{displaymath}
\lim_{N \rightarrow +\infty}\P \left( \left| \frac{X(\frac{[Nt]+1}{N},\frac{[Nt]+1}{N}) - X(\frac{[Nt]+1}{N},t_0)}{(1/N)^{H(t_0)}} \right| > u^{1/p}\right) =0.
\end{displaymath}
There exists $K_{U,p} >0$ such that 
\begin{equation}\label{majprobpart}
\P \left( \left| \frac{X(\frac{[Nt]+1}{N},\frac{[Nt]+1}{N}) - X(\frac{[Nt]+1}{N},t_0)}{(1/N)^{H(t_0)}} \right| > u^{1/p}\right) \leq \mathbf{1}_{u<1} + K_{U,p} \left(\frac{|\log u |^d}{u^{d/p}} + \frac{|\log u|^c}{u^{c/p}} \right) \mathbf{1}_{u \geq 1}.
\end{equation}

Since $\alpha$ is a continuous function, we can fix $U$ small enough such that $c=\inf_{t \in U}\limits \alpha(t) >p$. We deduce from the dominated convergence theorem that for all $t \in U$,
\begin{displaymath}
\lim_{N \rightarrow +\infty} \E \left[\left| \frac{X(\frac{[Nt]+1}{N},\frac{[Nt]+1}{N}) - X(\frac{[Nt]+1}{N},t_0)}{(1/N)^{H(t_0)}}\right|^p\right]  =0.
\end{displaymath}
With the inequality (\ref{majprobpart}),
\begin{displaymath}
\E \left[\left| \frac{X(\frac{[Nt]+1}{N},\frac{[Nt]+1}{N}) - X(\frac{[Nt]+1}{N},t_0)}{(1/N)^{H(t_0)}}\right|^p\right] \leq 1+\int_{1}^{+\infty}K_{U,p} \left(\frac{|\log u |^d}{u^{d/p}} + \frac{|\log u|^c}{u^{c/p}} \right) du
\end{displaymath}
and again with the dominated convergence theorem,
\begin{displaymath}
\lim_{N \rightarrow +\infty} \P \left(A_N(p) > x \right) =0.
\end{displaymath}
The same inequalities hold with $B_N(p)$ so we obtain  $B_N(p) \overset{\P}{\longrightarrow} 0$. We conclude proving  $C_N(p) \overset{\P}{\longrightarrow} \E|Z|^{p}$. Let $c_0 >0$. We use the decomposition

\begin{displaymath}
 C_N(p) - \E|Z|^{p} = \frac{1}{n(N)} \sum_{k=[Nt_0]-\frac{n(N)}{2}}^{[Nt_0]+\frac{n(N)}{2}-1} \left|\frac{X(\frac{k+1}{N},t_0)-X(\frac{k}{N},t_0)}{(1/N)^{H(t_0)}}\right|^p\mathbf{1}_{\left| \frac{X(\frac{k+1}{N},t_0)-X(\frac{k}{N},t_0)}{(1/N)^{H(t_0)}}\right|>c_0}- \E|Z|^{p}\mathbf{1}_{|Z|>c_0}
\end{displaymath}
\begin{displaymath}
 + \frac{1}{n(N)} \sum_{k=[Nt_0]-\frac{n(N)}{2}}^{[Nt_0]+\frac{n(N)}{2}-1} \left|\frac{X(\frac{k+1}{N},t_0)-X(\frac{k}{N},t_0)}{(1/N)^{H(t_0)}}\right|^p\mathbf{1}_{\left| \frac{X(\frac{k+1}{N},t_0)-X(\frac{k}{N},t_0)}{(1/N)^{H(t_0)}}\right|\leq c_0}- \E|Z|^{p}\mathbf{1}_{|Z| \leq c_0}.
\end{displaymath}
Let $\varepsilon >0$ and $x >0$. By Markov's inequality, we have
\begin{eqnarray*}
\P_1 & = &\P \left(\frac{1}{n(N)} \sum_{k=[Nt_0]-\frac{n(N)}{2}}^{[Nt_0]+\frac{n(N)}{2}-1} \left|\frac{X(\frac{k+1}{N},t_0)-X(\frac{k}{N},t_0)}{(1/N)^{H(t_0)}}\right|^p\mathbf{1}_{\left| \frac{X(\frac{k+1}{N},t_0)-X(\frac{k}{N},t_0)}{(1/N)^{H(t_0)}}\right|>c_0} > \frac{x}{4} \right)\\
 & \leq & \frac{4}{x n(N)} \sum_{k=[Nt_0]-\frac{n(N)}{2}}^{[Nt_0]+\frac{n(N)}{2}-1} \E \left[ \left|\frac{X(\frac{k+1}{N},t_0)-X(\frac{k}{N},t_0)}{(1/N)^{H(t_0)}}\right|^p\mathbf{1}_{\left| \frac{X(\frac{k+1}{N},t_0)-X(\frac{k}{N},t_0)}{(1/N)^{H(t_0)}}\right|>c_0}\right].\\
 \end{eqnarray*}
 Since $X(.,t_0)$ is $H(t_0)$-self-similar with stationary increments,
 \begin{displaymath}
 \P_1 \leq \frac{4}{x} \E \left[ \left|X(1,t_0)\right|^p\mathbf{1}_{\left|X(1,t_0)\right|>c_0}\right]
 \end{displaymath}
 and
 \begin{displaymath}
 \E|Z|^{p}\mathbf{1}_{|Z|\leq c_0} = \frac{1}{n(N)} \sum_{k=[Nt_0]-\frac{n(N)}{2}}^{[Nt_0]+\frac{n(N)}{2}-1} \E \left[ \left|\frac{X(\frac{k+1}{N},t_0)-X(\frac{k}{N},t_0)}{(1/N)^{H(t_0)}}\right|^p\mathbf{1}_{\left| \frac{X(\frac{k+1}{N},t_0)-X(\frac{k}{N},t_0)}{(1/N)^{H(t_0)}}\right| \leq c_0}\right].
\end{displaymath}
 
 We fix $c_0$ large enough such that for all $N \in \mathbb{N}$, $\P_1 \leq \frac{\varepsilon}{2}$ and $\E|Z|^{p}\mathbf{1}_{|Z|>c_0} < \frac{x}{4} $.
Writing $K(x)=|x|^p \mathbf{1}_{|x|\leq c_0}$ and $\Delta X_{k,t_0} = X(k+1,t_0) - X(k,t_0)$, using Chebyshev's inequality, we get
\begin{eqnarray*}
 \P \left( \left|C_N(p)- \E|Z|^{p} \right|> x \right) & \leq & \frac{\varepsilon}{2} + \frac{4}{x^2 n(N)^2} \sum_{k,j=[Nt_0]-\frac{n(N)}{2}}^{[Nt_0]+\frac{n(N)}{2}-1} Cov\left(K(\Delta X_{k,t_0}),K(\Delta X_{j,t_0}) \right)\\
 & \leq & \frac{\varepsilon}{2} + \frac{4}{x^2}\frac{ Var\left( K(\Delta X_{0,t_0})\right)}{ n(N)} + \frac{4}{x^2}\frac{1}{n(N)} \sum_{j=1}^{n(N)-1} Cov\left(K(\Delta X_{0,t_0}),K(\Delta X_{j,t_0}) \right).\\
\end{eqnarray*}
Under the condition {\bf (C*)}, we can apply Theorem 2.1 of \cite{PTA}: there exists a positive constant $C$ such that
\begin{displaymath}
|Cov\left(K(\Delta X_{0,t_0}),K(\Delta X_{j,t_0}) \right)| \leq C \Vert K \Vert_1^2 \int_E | h_{0,t_0}(v) h_{j,t_0}(v)|^{\frac{\alpha(t_0)}{2}} m(dv).
\end{displaymath}
Since the process $X(.,t_0)$ is $H(t_0)$-self-similar with stationary increments, the constant $C$ does not depend on $k$, $j$. We then obtain the existence of a positive constant $C_{p,c_0}$ depending on $p$, $c_0$ and $x$ such that
\begin{displaymath}
 \P \left( \left|C_N(p)- \E|Z|^{p} \right|> x \right) \leq \frac{\varepsilon}{2} + \frac{C_{p,c_0}}{n(N)}\int_E | h_{0,t_0}(v)|^{\alpha(t_0)} m(dv) + \frac{C_{p,c_0}}{n(N)} \sum_{j=1}^{n(N)-1} \int_E | h_{0,t_0}(v) h_{j,t_0}(v)|^{\frac{\alpha(t_0)}{2}} m(dv).
\end{displaymath}

Since $\lim_{N \rightarrow +\infty}\limits n(N)=+\infty$ and $\lim_{j \rightarrow +\infty}\limits \int_{E} | h_{0,t_0}(v) h_{j,t_0}(v)|^{\frac{\alpha(t_0)}{2}} m(dv) = 0$, we conclude with Cesaro's theorem that there exists $N_0 \in \mathbb{N}$ such that for all $N \geq N_0$,
\begin{displaymath}
 \frac{C_{p,c_0}}{n(N)}\int_E | h_{0,t_0}(v)|^{\alpha(t_0)} m(dv) + \frac{C_{p,c_0}}{n(N)} \sum_{j=1}^{n(N)-1} \int_E | h_{0,t_0}(v) h_{j,t_0}(v)|^{\frac{\alpha(t_0)}{2}} m(dv) \leq \frac{\varepsilon}{2}
\end{displaymath}
and
\begin{displaymath}
  \P \left( \left|C_N(p)- \E|Z|^{p} \right|> x \right) \leq \varepsilon \hspace*{0.5cm} \Box
\end{displaymath}

{\bf Proof of Lemma \ref{ProbaRapport}}

Let $\mu >0$ and $\lambda \in (0,1/e)$. Since $X(.,\frac{k}{N})$ is $H(\frac{k}{N})$-self-similar with stationary increments, $ N^{H(\frac{k}{N})}(X(\frac{k+1}{N},\frac{k}{N})- X(\frac{k}{N},\frac{k}{N}) )$ is distributed as the $\alpha(\frac{k}{N})$-stable variable $X(1,\frac{k}{N})$. We deduce that there exists $K_U >0$ such that 
$$\P \left( \frac{|X(\frac{k+1}{N},\frac{k}{N} ) - X(\frac{k}{N},\frac{k}{N} ) |}{(1/N)^{H(\frac{k}{N})}} \leq \mu \right) \leq K_U \mu.$$

Then 
\begin{displaymath}
 \P \left( \frac{|X(\frac{k+1}{N},\frac{k+1}{N}) - X(\frac{k+1}{N},\frac{k}{N}) |}{|X(\frac{k+1}{N},\frac{k}{N}) - X(\frac{k}{N},\frac{k}{N}) |} > \lambda \right) \leq \P \left( \frac{|X(\frac{k+1}{N},\frac{k+1}{N}) - X(\frac{k+1}{N},\frac{k}{N}) |}{(1/N)^{H(\frac{k}{N})}} > \lambda \mu \right) + K_U \mu.
\end{displaymath}

With the conditions {\bf (R1)}, {\bf (M1)}, {\bf (M2)} and {\bf (M3)}, we use the inequality (\ref{majprob}) to obtain (with $K_U$ which may change from line to line)
\begin{displaymath}
 \P \left( \frac{|X(\frac{k+1}{N},\frac{k+1}{N}) - X(\frac{k+1}{N},\frac{k}{N}) |}{|X(\frac{k+1}{N},\frac{k}{N}) - X(\frac{k}{N},\frac{k}{N}) |} > \lambda \right) \leq K_U \frac{N^{dH(\frac{k}{N})}}{|N \lambda \mu|^d}(1+|\log |N \lambda \mu | |^d ) + K_U \frac{N^{cH(\frac{k}{N})}}{|N \lambda \mu|^c}(1+|\log |N \lambda \mu | |^c ) + K_U \mu.
\end{displaymath}
We choose $\mu = \frac{1}{\lambda^{\frac{\alpha(t_0)}{1+\alpha(t_0)}} N^{\frac{\alpha(t_0)(1-H(\frac{k}{N}))}{1+\alpha(t_0)}}}$ to obtain

\begin{eqnarray*}
 \P \left( \frac{|X(\frac{k+1}{N},\frac{k+1}{N}) - X(\frac{k+1}{N},\frac{k}{N}) |}{|X(\frac{k+1}{N},\frac{k}{N}) - X(\frac{k}{N},\frac{k}{N}) |} > \lambda \right) & \leq & K_U (\frac{|\log N|^d|\log \lambda |^d}{N^{\frac{d(1-H(\frac{k}{N}))}{1+\alpha(t_0)}}\lambda^{\frac{d}{1+\alpha(t_0)}}} 
 + \frac{|\log N|^c|\log \lambda |^c}{N^{\frac{c(1-H(\frac{k}{N}))}{1+\alpha(t_0)}} \lambda^{\frac{c}{1+\alpha(t_0)}}}) + K_U \mu\\
 & \leq & K_U \frac{|\log N|^d|\log \lambda |^d}{N^{\frac{d(1-H_{-})}{1+c}}\lambda^{\frac{d}{1+c}}} \Box \\ 
\end{eqnarray*}

{\bf Proof of Theorem \ref{ConvReste}}

Let $x >0$ and $\delta \in ( 0,\frac{2\alpha(t_0)(1-H(t_0))}{2+3\alpha(t_0)} )$. Put $\xi_{k,N} = \frac{X(\frac{k+1}{N},\frac{k+1}{N})-X(\frac{k+1}{N},\frac{k}{N})}{X(\frac{k+1}{N},\frac{k}{N})-X(\frac{k}{N},\frac{k}{N})}.$

Let us show that $\frac{1}{\sqrt{n(N)}} \sum_{k=[Nt_0]-\frac{n(N)}{2}}^{[Nt_0]+\frac{n(N)}{2}-1}\limits \log |1+\xi_{k,N}|$ tends to $0$ in probability. Since $\alpha$ and $H$ are continuous, we can choose $\mu > \frac{1}{2}$ and $U$ small enough in order to have $ \delta < \frac{d(1-H_{-})}{1+c+\mu d }.$

Let $\lambda_N=1-\frac{1}{n(N)^{\mu}}$, $\mu_N = \frac{1}{n(N)^{\mu}}$, $A_N= \cup_{k=[Nt_0]-\frac{n(N)}{2}}^{[Nt_0]+ \frac{n(N)}{2}-1} \{ |1+\xi_{k,N} | < \lambda_N \}$ and $B_N = \cup_{k=[Nt_0]-\frac{n(N)}{2}}^{[Nt_0]+ \frac{n(N)}{2}-1} \{ | \xi_{k,N} | > \mu_N \}$.
Since $A_N \subset B_N$, we will only show that $\P(B_N) =0$. We use Lemma \ref{ProbaRapport}: there exists $K_U >0$ such that 
\begin{eqnarray*}
 \P(|\xi_{k,N}| > \mu_N) & \leq & K_U \frac{|\log N|^d n(N)^{\frac{d \mu}{1+c}}}{N^{\frac{d(1-H_{-})}{1+c}}}. \\
 \end{eqnarray*}
 Then
 \begin{eqnarray*}
\P(B_N) & \leq &  n(N) K_U \frac{|\log N|^d n(N)^{\frac{d \mu}{1+c}}}{N^{\frac{d(1-H_{-})}{1+c}}} \\
& \leq & K_U |\log N|^d N^{\delta (1+ \frac{d \mu}{1+c}) - \frac{d(1-H_{-})}{1+c}} \\
\end{eqnarray*}
so $\lim_{N \rightarrow +\infty}\limits \P(B_N) =0.$ 

We obtain then

\begin{eqnarray*}
 \P \left(\sum_{k=[Nt_0]-\frac{n(N)}{2}}^{[Nt_0]+\frac{n(N)}{2}-1}\limits \frac{\log |1+\xi_{k,N}|}{\sqrt{n(N)}} < - x \right) & \leq & \P(A_N) + \P \left( \{ \sum_{k=[Nt_0]-\frac{n(N)}{2}}^{[Nt_0]+\frac{n(N)}{2}-1}\limits \frac{\log |1+\xi_{k,N}|}{\sqrt{n(N)}} < - x \} \cap  \bar{A_N}\right) \\
& \leq  & \P(A_N) + \P \left( n(N) \frac{\log \lambda_N}{\sqrt{n(N)}} < -x\right).\\
 \end{eqnarray*}
 Since $\mu > \frac{1}{2}$, $\lim_{N \rightarrow +\infty}\limits \P \left( \sqrt{n(N)} \log \lambda_N < -x\right)=0.$ We obtain in the same way
 \begin{eqnarray*}
 \P \left(\sum_{k=[Nt_0]-\frac{n(N)}{2}}^{[Nt_0]+\frac{n(N)}{2}-1}\limits \frac{\log |1+\xi_{k,N}|}{\sqrt{n(N)}} > x \right) & \leq & \P(B_N) + \P \left( \{ \sum_{k=[Nt_0]-\frac{n(N)}{2}}^{[Nt_0]+\frac{n(N)}{2}-1}\limits \frac{\log |1+\xi_{k,N}|}{\sqrt{n(N)}} > x \} \cap  \bar{B_N}\right) \\
& \leq  & \P(B_N) + \P \left( n(N) \frac{\log |1+ \mu_N|}{\sqrt{n(N)}} >x\right).\\
 \end{eqnarray*}
 Since $\mu > \frac{1}{2}$, $\lim_{N \rightarrow +\infty}\limits \P \left( \sqrt{n(N)} \log |1+\mu_N| >x\right)=0$ \Box

 {\bf Proof of Theorem \ref{ConvLpalpha}}

Since $x \rightarrow x^{\gamma}$ is an increasing function on $\bbbr_+$ (we take $\gamma \in (0,1)$), 
\begin{displaymath}
 \hat{\alpha}_N(t_0) = \min \left(\argmin_{\alpha \in [0,2]} \int_{p_0}^{2} |R_{\textrm{exp}}^{(N)}(p) - R_{\alpha}(p) |^{\gamma} dp \right).
\end{displaymath}
Let $g_N(\alpha) = \int_{p_0}^{2}\limits |R_{\textrm{exp}}^{(N)}(p) - R_{\alpha}(p) |^{\gamma} dp $ and $g(\alpha) =  \int_{p_0}^{2}\limits |R_{\alpha(t_0)}(p) - R_{\alpha}(p) |^{\gamma} dp$. 



$g$ is a continuous function on $(0,2]$, with $g(0)>0$, $g(2)>0$. The only solution of the equation $g(\alpha)=0$ is $\alpha(t_0)$. Moreover, $\lim_{\alpha \rightarrow \alpha(t_0)}\limits \frac{|g(\alpha) - g(\alpha(t_0) |}{|\alpha - \alpha(t_0)|^{\gamma}} >0$.



Then, there exists  $K_{\alpha(t_0)}$  a positive constant depending only on $\alpha(t_0)$ such that:
\begin{equation}\label{inegalpha}
 \forall \alpha \in (0,2), \hspace{0.2cm} |g(\alpha)| \geq K_{\alpha(t_0)} |\alpha- \alpha(t_0)|.
\end{equation}
 We estimate now $|g(\hat{\alpha}_N(t_0))|$.
\begin{eqnarray*}
 |g(\hat{\alpha}_N(t_0))| & \leq & |g(\hat{\alpha}_N(t_0) ) - g_N(\hat{\alpha}_N(t_0)) | + |g_N(\hat{\alpha}_N(t_0))|\\
 & \leq & |g(\hat{\alpha}_N(t_0) ) - g_N(\hat{\alpha}_N(t_0)) | + g_N(\alpha(t_0)),\\
\end{eqnarray*}
and
\begin{eqnarray*}
 |g(\hat{\alpha}_N(t_0) ) - g_N(\hat{\alpha}_N(t_0)) | & = & \left| \int_{p_0}^2 \left( |R_{\alpha(t_0)}(p)-R_{\hat{\alpha}_N(t_0)}(p) |^{\gamma} - |R_{\textrm{exp}}^{(N)}(p) - R_{\hat{\alpha}_N(t_0)}(p) |^{\gamma}\right) dp \right|\\
 & \leq & \int_{p_0}^2 \left|R_{\alpha(t_0)}(p)-R_{\textrm{exp}}^{(N)}(p) \right|^{\gamma} dp \\
 & = & g_N(\alpha(t_0)).\\ 
\end{eqnarray*}
From (\ref{inegalpha}),
\begin{eqnarray*}
 \left| \hat{\alpha}_N(t_0) - \alpha(t_0) \right| & \leq & \frac{1}{K_{\alpha(t_0)}}g(\hat{\alpha}_N(t_0) ) \\
 & \leq & \frac{2}{K_{\alpha(t_0)}} g_N(\alpha(t_0)).
\end{eqnarray*}

Let us show that $\lim_{N \rightarrow +\infty}\limits \E \left|g_N(\alpha(t_0))\right|^r = 0$ for any $r >0$. One has, using the inequality $S_N(p) \leq S_N(q)$ for $p\leq q$,
\begin{eqnarray*}
  g_N(\alpha(t_0)) & = & \int_{p_0}^{\alpha(t_0)}\limits |R_{\textrm{exp}}^{(N)}(p) - R_{\alpha(t_0)}(p) |^{\gamma} dp + \int_{\alpha(t_0)}^{2}\limits |R_{\textrm{exp}}^{(N)}(p)|^{\gamma} dp\\
  & \leq & \int_{p_0}^{\alpha(t_0)}\limits |R_{\textrm{exp}}^{(N)}(p) - R_{\alpha(t_0)}(p) |^{\gamma} dp + (2-\alpha(t_0))\left| \frac{S_N(p_0)}{S_N(\alpha(t_0))}\right|^{\gamma}.\\
\end{eqnarray*}

For the first term, we use Theorem \ref{ConvLpSnp} : 
for all $p \in [p_0,\alpha(t_0))$,
\begin{equation}\label{convprob}
  N^{H(t_0)}S_N(p) \overset{\P}{\longrightarrow} (\E|X(1,t_0)|^{p})^{1/p}
\end{equation}
It is clear that $\forall p \in [p_0,\alpha(t_0))$,
\begin{displaymath}
 \left(N^{H(t_0)}S_N(p_0),N^{H(t_0)}S_N(p)\right) \overset{\P}{\longrightarrow} \left( (\E|X(1,t_0)|^{p_0})^{1/p_0}, (\E|X(1,t_0)|^{p})^{1/p}\right),
\end{displaymath}
and
\begin{equation}\label{convprobaRexp}
 R_{\textrm{exp}}^{(N)}(p) = \frac{S_N(p_0)}{S_N(p)} \overset{\P}{\longrightarrow} R_{\alpha(t_0)}(p).
\end{equation}

Note that $\forall N \in \mathbb{N}$, $\forall p \in [p_0,\alpha(t_0))$, $|R_{\textrm{exp}}^{(N)}(p)| \leq 1$ so there exists a positive constant $K$ depending on $\gamma r$, $\alpha(t_0)$ and $p$ such that
\begin{eqnarray*}
\E |R_{\textrm{exp}}^{(N)}(p) - R_{\alpha(t_0)}(p) |^{\gamma r} & = & \int_{0}^{K} \P \left( |R_{\textrm{exp}}^{(N)}(p) - R_{\alpha(t_0)}(p) |^{\gamma r} >x \right) dx.\\
\end{eqnarray*}
Finally, with (\ref{convprobaRexp}),  $\forall p \in [p_0,\alpha(t_0))$, $\E |R_{\textrm{exp}}^{(N)}(p) - R_{\alpha(t_0)}(p) |^{\gamma r} \underset{N \rightarrow +\infty}{\longrightarrow} 0.$ With the inequality $\E |R_{\textrm{exp}}^{(N)}(p) - R_{\alpha(t_0)}(p) |^{\gamma r} \leq 2C_{\gamma r}$ where $C_{\gamma r}$ is a positive constant depending on $\gamma r$, by the dominating convergence theorem,
\begin{displaymath}
 \lim_{N \rightarrow +\infty} \int_{p_0}^{\alpha(t_0)} \E |R_{\textrm{exp}}^{(N)}(p) - R_{\alpha(t_0)}(p) |^{\gamma r} dp =0.
\end{displaymath}
To conclude we show that $\left| \frac{S_N(p_0)}{S_N(\alpha(t_0))}\right|^{\gamma}  \overset{L^r}{\longrightarrow} 0$. Since $\forall N \in \mathbb{N}$, $\left| \frac{S_N(p_0)}{S_N(\alpha(t_0))}\right|^{\gamma} \leq 1$, it is enough to show $\frac{S_N(p_0)}{S_N(\alpha(t_0))} \overset{\P}{\longrightarrow} 0$. Let $p < \alpha(t_0)$.
\begin{displaymath}
 \P(\frac{1}{N^{H(t_0)}S_N(\alpha(t_0)) } > x) \leq \P(\frac{1}{N^{H(t_0)}S_N(p) } > x).
\end{displaymath}

So,
\begin{eqnarray*}
 \limsup_{N \rightarrow +\infty} \P(\frac{1}{N^{H(t_0)}S_N(\alpha(t_0))} > x) & \leq & \limsup_{N \rightarrow +\infty} \P(\frac{1}{N^{H(t_0)}S_N(p)} > x)\\
 & = & \lim_{N \rightarrow +\infty} \P(\frac{1}{N^{H(t_0)}S_N(p) } > x)\\
 & = & \P(\frac{1}{ (\E|X(1,t_0)|^{p})^{1/p}} > x),\\
\end{eqnarray*}
with (\ref{convprob}). Since $\lim_{p \rightarrow \alpha(t_0)}\limits \P(\frac{1}{ (\E|X(1,t_0)|^{p})^{1/p}} > x) =0$, we have  $\limsup_{N \rightarrow +\infty}\limits \P(\frac{1}{N^{H(t_0)}S_N(\alpha(t_0))} > x ) =0$ and $\frac{1}{N^{H(t_0)}S_N(\alpha(t_0))} \overset{\P}{\longrightarrow} 0$.
Using the convergence  $N^{H(t_0)}S_N(p_0) \overset{\P}{\longrightarrow} (\E|X(1,t_0)|^{p_0})^{1/p_0} $, we obtain $\frac{S_N(p_0)}{S_N(\alpha(t_0))} \overset{\P}{\longrightarrow} 0$.

If in addition, we assume that the process $X(.,t_0)$ is $H(t_0)$-self-similar with stationary increments and $H(t_0) <1$, and $\lim_{j \rightarrow +\infty}\limits \int_E | h_{0,t_0}(x) h_{j,t_0}(x)|^{\frac{\alpha(t_0)}{2}} m(dx)=0$, for all $t_0 \in U$, we obtain for all $r>0$ and all $t \in U$
$$\lim_{N \rightarrow +\infty}\E \left|\hat{\alpha}_N(t)- \alpha(t) \right|^r = 0.$$
$\hat{\alpha}_N$ and $\alpha$ are two bounded functions on $U$ so for all $r>0$, 
$$\lim_{N \rightarrow +\infty}\limits \int_U\limits \E \left[\left|\hat{\alpha}_N(t)- \alpha(t) \right|^r\right] dt = 0 \quad \Box$$

{\bf Proof of Theorem \ref{ConvLpH}}

Note that it is sufficient to prove the result of Theorem \ref{ConvLpH} for $r \geq 1$ since the convergence in $L^p$ implies the convergence in $L^q$ for all $q<p$. Let $r \geq 1$.
We write 
\begin{eqnarray*}
 \hat{H}_N(t_0) - H(t_0) & = & -\frac{1}{n(N) \log N} \sum_{k=[Nt_0]-\frac{n(N)}{2}}^{[Nt_0]+\frac{n(N)}{2}-1} \log \left|\frac{Y_{k,N}}{(\frac{1}{N})^{H(t_0)}} \right|\\
 & =& -\frac{N}{n(N) \log N} \int_{\frac{[Nt_0]}{N}-\frac{n(N)}{2N}}^{\frac{[Nt_0]}{N}+\frac{n(N)}{2N} } \log \left| \frac{Y(\frac{[Nt]+1}{N}) - Y(\frac{[Nt]}{N})}{(\frac{1}{N})^{H(t_0)}}\right|  dt.\\
\end{eqnarray*}
Let $\delta_N(dt) = \frac{N}{n(N)} \mathbf{1}_{\{ \frac{[Nt_0]}{N}-\frac{n(N)}{2N} \leq t < \frac{[Nt_0]}{N}+\frac{n(N)}{2N}\}} dt $ and $ f_N(t) = \log \left| \frac{Y(\frac{[Nt]+1}{N}) - Y(\frac{[Nt]}{N})}{(\frac{1}{N})^{H(t)}}\right|$.\par Since $\int_{0}^{1} \delta_N(dt) =1$,  we obtain
\begin{displaymath}
 \hat{H}_N(t_0) - H(t_0) = -\frac{1}{\log N} \int_{0}^{1} f_N(t) \delta_N(dt) + \int_{0}^{1} \left( H(t)-H(t_0)\right) \delta_N(dt). 
\end{displaymath}
Then, there exists a constant $K_r >0$ depending on $r$ such that
\begin{displaymath}
 \E \left[ |\hat{H}_N(t_0) - H(t_0)|^r \right] \leq  K_r \frac{\E \left( |\int_{0}^{1} f_N(t) \delta_N(dt) |^r\right)}{|\log N|^r}  + K_r \left| \int_{0}^{1} \left( H(t)-H(t_0)\right) \delta_N(dt) \right|^r.
\end{displaymath}

$H$ is continuously differentiable and $\lim_{N \rightarrow +\infty}\limits \frac{N}{n(N)} = + \infty$ so 
\begin{displaymath}
 \lim_{N \rightarrow + \infty} \int_{0}^{1} \left( H(t)-H(t_0)\right) \delta_N(dt) = 0.
\end{displaymath}
To conclude, it is sufficient to show that there exists a constant $K >0$ depending on $t_0$ and $r$ such that for all $N \in \mathbb{N}$, $\E \left( |\int_{0}^{1} f_N(t) \delta_N(dt) |^r\right) \leq K$. 
Let $U$ an open interval satisfying all the conditions {\bf (R-)}, {\bf (M-)} and {\bf (H-)}, and $t_0 \in U$. We can fix $N_0 \in \mathbb{N}$ and $V \subset U$ an open interval depending on $t_0$ such that for all $N \geq N_0$ and all $t \in V$, $ \frac{[Nt]+1}{N} \in U$, $\frac{[Nt]}{N} \in U$ and $\int_{0}^{1} f_N(t) \delta_N(dt) =  \int_V f_N(t) \delta_N(dt)$.
 With the Jensen inequality,
 \begin{displaymath}
 \E \left( |\int_{0}^{1} f_N(t) \delta_N(dt) |^r\right) \leq \int_V \E |f_N(t)|^r \delta_N(dt).
\end{displaymath}
We consider $\E |f_N(t)|^r = \int_{0}^{+\infty}\limits \P \left(| f_N(t)|^r > x \right)dx$.

\begin{eqnarray*}
 \E |f_N(t)|^r & = & \int_{0}^{+\infty}\limits \P \left( \left| Y(\frac{[Nt]+1}{N}) - Y(\frac{[Nt]}{N})\right| < \frac{e^{-x^{1/r}}}{N^{H(t)}} \right)dx.\\
 & & +\int_{0}^{+\infty}\limits \P \left( \left| Y(\frac{[Nt]+1}{N}) - Y(\frac{[Nt]}{N})\right| > \frac{e^{x^{1/r}}}{N^{H(t)}} \right)dx \\
\end{eqnarray*}

Thanks to the conditions {\bf (R1)}, {\bf (M4)}, {\bf (M5)}, {\bf (M6)}, {\bf (M7)}, {\bf (H1)}, {\bf (H3)}, {\bf (H4)} and {\bf (H5)}, we use the equality (\ref{probinferi}) to control the first term: there exists $K_U >0$ (that may change from line to line) and $N_0 \in \bbbn$ such that for all $N \geq N_0$ and all $t \in V$,

\begin{displaymath}
 \P \left( | Y(\frac{[Nt]}{N}+\frac{1}{N}) - Y(\frac{[Nt]}{N}) | < \frac{e^{-x^{1/r}}}{N^{H(t)}}\right) \leq K_U N^{H(\frac{[Nt]}{N})} \frac{e^{-x^{1/r}}}{N^{H(t)}}.
\end{displaymath}
We get then
\begin{eqnarray*}
 \int_{0}^{+\infty}\limits \P \left( \left| Y(\frac{[Nt]+1}{N}) - Y(\frac{[Nt]}{N})\right| < \frac{e^{-x^{1/r}}}{N^{H(t)}} \right)dx & \leq & K_U (\int_{0}^{+\infty}\limits e^{-x^{1/r}} dx  ) N^{H(\frac{[Nt]}{N})-H(t)}\\
 & \leq & K_U.\\
\end{eqnarray*}

For the second term, we write
\begin{eqnarray*}
 \P \left( \left| Y(\frac{[Nt]+1}{N}) - Y(\frac{[Nt]}{N})\right| > \frac{e^{x^{1/r}}}{N^{H(t)}} \right) & \leq& \P \left( \left| X(\frac{[Nt]+1}{N},\frac{[Nt]+1}{N}) - X(\frac{[Nt]+1}{N},\frac{[Nt]}{N})\right| > \frac{e^{x^{1/r}}}{N^{H(t)}} \right)\\
 & & +\P \left( \left| X(\frac{[Nt]+1}{N},\frac{[Nt]}{N}) - X(\frac{[Nt]}{N},\frac{[Nt]}{N})\right| > \frac{e^{x^{1/r}}}{N^{H(t)}} \right) .\\
\end{eqnarray*}
With the conditions {\bf (R1)}, {\bf (M1)}, {\bf (M2)} and {\bf (M3)}, we use the equality (\ref{majprob}) to obtain a positive constant $K_U >0$ such that:

\begin{eqnarray*}
& & \P \left( \left| X(\frac{[Nt]+1}{N},\frac{[Nt]+1}{N}) - X(\frac{[Nt]+1}{N},\frac{[Nt]}{N})\right| > \frac{e^{x^{1/r}}}{N^{H(t)}} \right)\\
& & \leq  K_U \left( \frac{(\log N)^c}{N^{c(1-H(t))}e^{cx^{1/r}}}\right) + K_U \left( \frac{x^{c/r}}{N^{c(1-H(t))}e^{cx^{1/r}}}\right)\\
 & &   + K_U \left(\frac{(\log N)^d}{N^{d(1-H(t))}e^{dx^{1/r}}} \right) + K_U \left( \frac{x^{d/r}}{N^{d(1-H(t))}e^{dx^{1/r}}}\right).\\
\end{eqnarray*}
Since $H_+ <1$, we conclude that
\begin{displaymath}
 \lim_{N \rightarrow +\infty}  \int_{0}^{+\infty} \limits \sup_{t \in U}\limits \P \left( \left| X(\frac{[Nt]+1}{N},\frac{[Nt]+1}{N} ) - X(\frac{[Nt]+1}{N}, \frac{[Nt]}{N})\right| \geq \frac{e^{x^{1/r}}}{N^{H(t)}} \right) dx = 0.
\end{displaymath}
Let $\eta < c$. The Markov inequality gives
\begin{displaymath}
 \P \left( \left| X(\frac{[Nt]+1}{N},\frac{[Nt]}{N} ) - X(\frac{[Nt]}{N}, \frac{[Nt]}{N})\right| \geq \frac{e^{x^{1/r}}}{N^{H(t)}} \right) \leq \frac{N^{\eta H(t)}}{e^{\eta x^{1/r}}} \E \left[|X(\frac{[Nt]+1}{N},\frac{[Nt]}{N}) -X(\frac{[Nt]}{N},\frac{[Nt]}{N}) |^{\eta}\right]
\end{displaymath}
and Property 1.2.17 of \cite{ST}
\begin{displaymath}
\E \left[|X(\frac{[Nt]+1}{N},t_N) -X(t_N,t_N) |^{\eta}\right] = c_{\alpha(t_N),0}(\eta)^{\eta} \left(\int_E (|f(\frac{[Nt]+1}{N},t_N,x) - f(t_N,t_N,x)|^{\alpha(t_N)} \hspace{0.1cm}  m(dx) \right)^{\eta/\alpha(t_N)}.
\end{displaymath}
where $t_N=\frac{[Nt]}{N}$. With the condition {\bf (H2)}, there exists $K_U > 0$ such that for all $N \geq N_0$ and all $t \in V$,

\begin{displaymath}
 \int_{0}^{+\infty} \limits \P \left( \left| X(\frac{[Nt]+1}{N},\frac{[Nt]}{N} ) - X(\frac{[Nt]}{N},\frac{[Nt]}{N})\right| \geq \frac{e^{x^{1/r}}}{N^{H(t)}} \right) dx \leq K_U.
\end{displaymath}

The conclusion is that for all $t_0 \in U$, 
 \begin{equation}\label{limennorme}
 \lim_{N \rightarrow +\infty}\limits \E \left|\hat{H}_N(t_0)- H(t_0) \right|^r = 0.
\end{equation}
Let $[a,b] \subset U$, $p>0$ and $\eta \in (0,1)$. We denote $A=\liminf_{N \rightarrow + \infty}\limits \{ \sup_{t \in [a,b]}\limits | \hat{H}_N(t) - H(t)| \leq B\}.$ Thanks to Lemma \ref{LemSup}, there exists $B \in \bbbr$ such that $\P (A)=1$. Then 
\begin{eqnarray*}
 \E \left[ \int_a^b\limits |\hat{H}_N(t)- H(t) |^p dt\right] & = &  \int_a^b\E \left[  |\hat{H}_N(t)- H(t) |^p \right] dt\\
 & = & \int_a^b \int_{0}^{+\infty} \P \left( \{|\hat{H}_N(t)- H(t) |^p > x \} \cap A \right) dx \hspace{0.1cm} dt\\
 & = & \int_a^b \int_{0}^{B^p} \P \left( \{|\hat{H}_N(t)- H(t) |^p > x \} \cap A \right) dx \hspace{0.1cm} dt \\
 & \leq & \int_a^b \int_{0}^{B^p} \P \left( \{|\hat{H}_N(t)- H(t) |^p > x \} \right) dx \hspace{0.1cm} dt \\
\end{eqnarray*}
The equality (\ref{limennorme}) available for all $r >0$ easily leads to 
\begin{displaymath}
 \lim_{N \rightarrow +\infty}\limits \E \left[ \int_a^b\limits |\hat{H}_N(t)- H(t) |^p dt\right] = 0 \quad \Box
\end{displaymath}

{\bf Proof of Theorem \ref{Convlaw}}

We write $H_N^s(t_0) = -\frac{1}{n(N) \log N} \sum_{k=[Nt_0] - \frac{n(N)}{2}}^{[Nt_0] + \frac{n(N)}{2}-1}\limits \log \left| X(\frac{k+1}{N},\frac{k}{N})-X(\frac{k}{N},\frac{k}{N})\right|$ and the following decomposition
$$\log N (\hat{H}_N(t_0) - H(t_0)) + \mu_{t_0} = \log N (\hat{H}_N(t_0) - H_N^s(t_0)) + (\log N (H_N^s(t_0) - H(t_0)) + \mu_{t_0} ).$$
We know from \cite{LGLV3} that $Y$ is satisfying the conditions ${\bf (R1)}$, ${\bf (M1)}$, ${\bf (M2)}$ and ${\bf (M3)}$. For all $u \in (0,1)$, $X(.,u)$ is a $\alpha(u)$-stable L\'evy motion, so is $\frac{1}{\alpha(u)}$-self-similar with stationary increments. 
We apply Theorem \ref{ConvReste} to obtain the convergence in probability to $0$ of $\sqrt{n(N)} \log N (\hat{H}_N(t_0) - H_N^s(t_0) )$. For the second term, notice that

\begin{eqnarray*}
 H_N^s(t_0) & = & -\frac{1}{n(N) \log N} \sum_{k=[Nt_0] - \frac{n(N)}{2}}^{[Nt_0] + \frac{n(N)}{2}-1}\limits \log \left| X(\frac{k+1}{N},\frac{k}{N})-X(\frac{k}{N},\frac{k}{N})\right|\\
&= & -\frac{1}{n(N) \log N} \sum_{k=[Nt_0] - \frac{n(N)}{2}}^{[Nt_0] + \frac{n(N)}{2}-1}\limits \log \left| \frac{ X(\frac{k+1}{N},\frac{k}{N})-X(\frac{k}{N},\frac{k}{N})}{(1/N )^{H(k/N)}}\right| + \frac{1}{n(N)} \sum_{k=[Nt_0] - \frac{n(N)}{2}}^{[Nt_0] + \frac{n(N)}{2}-1}\limits H(\frac{k}{N}).\\
\end{eqnarray*}

Put $z_{k,N} = \log \left| \frac{ X(\frac{k+1}{N},\frac{k}{N})-X(\frac{k}{N},\frac{k}{N})}{(1/N )^{H(k/N)}}\right|$. Then 
\begin{eqnarray*}
  \sqrt{n(N)}(\log N (H_N^s(t_0) - H(t_0)) + \mu_{t_0} )  = && \frac{\log N}{\sqrt{n(N)}} \sum_{k=[Nt_0] - \frac{n(N)}{2}}^{[Nt_0] + \frac{n(N)}{2}-1}\limits (H(\frac{k}{N}) - H(t_0))\\
& & + \frac{1}{\sqrt{n(N)}} \sum_{k=[Nt_0] - \frac{n(N)}{2}}^{[Nt_0] + \frac{n(N)}{2}-1}\limits (\mu_{t_0} - \mu_{\frac{k}{N}} )\\
& &+ \frac{1}{\sqrt{n(N)}} \sum_{k=[Nt_0] - \frac{n(N)}{2}}^{[Nt_0] + \frac{n(N)}{2}-1}\limits(\mu_{\frac{k}{N}}-z_{k,N} ).\\
 \end{eqnarray*}
$H = \frac{1}{\alpha}$ is a $\mathcal{C}^1$ function, so there exists $K>0$ such that $|H(\frac{k}{N}) - H(t_0)| \leq K|\frac{k}{N} - t_0 |$, and

\begin{eqnarray*}
 \left| \frac{\log N}{\sqrt{n(N)}} \sum_{k=[Nt_0] - \frac{n(N)}{2}}^{[Nt_0] + \frac{n(N)}{2}-1}\limits (H(\frac{k}{N}) - H(t_0)) \right| & \leq & K \frac{\log N}{\sqrt{n(N)}} \sum_{k=[Nt_0] - \frac{n(N)}{2}}^{[Nt_0] + \frac{n(N)}{2}-1}\limits |\frac{k}{N} - t_0 | \\
 & \leq & K \frac{n(N)^{\frac{3}{2}}}{N} \log N.\\
\end{eqnarray*}
Since $\delta < \frac{2 \alpha(t_0) -2}{3 \alpha(t_0) +2} $, $\frac{3}{2}\delta <1$ and $\lim_{N \rightarrow +\infty}\limits \frac{\log N}{\sqrt{n(N)}} \sum_{k=[Nt_0] - \frac{n(N)}{2}}^{[Nt_0] + \frac{n(N)}{2}-1}\limits (H(\frac{k}{N}) - H(t_0)) =0.$

With $Z \sim S_{\alpha}(1,0,0)$, we can use the inversion formula to obtain the equality $$\E [ \log |Z| ] = \int_{\bbbr} \log |x| \frac{1}{\pi} \int_{0}^{+\infty} e^{-|t|^{\alpha}} \cos(tx) \hspace{0.1cm} dt \hspace{0.1cm} dx,$$ and check that the function $\alpha \mapsto \E [ \log |Z| ]$ is continuously differentiable.
 With the hypothesis on the function $\alpha$, the function $t \mapsto \mu_t$ is a $\mathcal{C}^1$ function. We get then, as for $H$, $ \lim_{N \rightarrow +\infty}\limits \frac{1}{\sqrt{n(N)}} \sum_{k=[Nt_0] - \frac{n(N)}{2}}^{[Nt_0] + \frac{n(N)}{2}-1}\limits (\mu_{\frac{k}{N}} - \mu_{t_0}) =0.$
To finish the proof, let us show the convergence $ \frac{1}{\sqrt{n(N)}} \sum_{k=[Nt_0] - \frac{n(N)}{2}}^{[Nt_0] + \frac{n(N)}{2}-1}\limits(\mu_{\frac{k}{N}}-z_{k,N} ) \cd \mathcal{N}(0,\sigma^2_{t_0}).$
Let $X_{k,N} = \frac{\mu_{\frac{k}{N}}-z_{k,N}}{\sqrt{n(N)}}$, $\varepsilon >0$ and $c = \inf_{t \in U} \alpha(t)$.

\begin{eqnarray*}
 \P \left( |X_{k,N}| > \varepsilon \right) & = & \P \left( \left| \frac{ X(\frac{k+1}{N},\frac{k}{N})-X(\frac{k}{N},\frac{k}{N})}{(1/N )^{H(k/N)}}\right| \leq  e^{\mu_{\frac{k}{N}}} e^{-\varepsilon \sqrt{n(N)}} \right)  \\
&  & + \P \left( \left| \frac{ X(\frac{k+1}{N},\frac{k}{N})-X(\frac{k}{N},\frac{k}{N})}{(1/N )^{H(k/N)}}\right|^{\frac{c}{2}} > e^{\frac{c }{2}\mu_{\frac{k}{N}}} e^{\frac{c\varepsilon \sqrt{n(N)}}{2}} \right). \\
 \end{eqnarray*}

$\frac{ X(\frac{k+1}{N},\frac{k}{N})-X(\frac{k}{N},\frac{k}{N})}{(1/N )^{H(k/N)}}$ is a standard $\alpha(\frac{k}{N})$-stable random variable, then there exists $K>0$ such that 
\begin{equation}\label{probaepsilon}
 \P \left( |X_{k,N}| > \varepsilon \right) \leq K ( e^{-\frac{c\varepsilon \sqrt{n(N)}}{2}} + e^{-\varepsilon \sqrt{n(N)}} ).
\end{equation}

$(X_{k,N})_k$ is thus satisfying $\lim_{N \rightarrow +\infty}\limits \sum_{k=[Nt_0] - \frac{n(N)}{2}}^{[Nt_0] + \frac{n(N)}{2}-1}\limits \P(|X_{k,N}| > \varepsilon ) =0.$ 

$\mu_{\frac{k}{N}} = \E [z_{k,N} ]$ so $$\sum_{k=[Nt_0] - \frac{n(N)}{2}}^{[Nt_0] + \frac{n(N)}{2}-1}\limits \E [X_{k,N} \mathbf{1}_{|X_{k,N}| \leq \varepsilon} ] = - \sum_{k=[Nt_0] - \frac{n(N)}{2}}^{[Nt_0] + \frac{n(N)}{2}-1}\limits \E [ X_{k,N} \mathbf{1}_{|X_{k,N}|> \varepsilon}].$$
With the inequality (\ref{probaepsilon}), we obtain $\lim_{N \rightarrow +\infty}\limits  \sum_{k=[Nt_0] - \frac{n(N)}{2}}^{[Nt_0] + \frac{n(N)}{2}-1}\limits \E [X_{k,N} \mathbf{1}_{|X_{k,N}| \leq \varepsilon} ] =0.$

Finally, we have $$\sum_{k=[Nt_0] - \frac{n(N)}{2}}^{[Nt_0] + \frac{n(N)}{2}-1}\limits \Var(X_{k,N} \mathbf{1}_{|X_{k,N}| \leq \varepsilon}) = \int_0^1 \Var \left((\mu_{\frac{[Nt]}{N}}-z_{[Nt],N})\mathbf{1}_{|X_{[Nt],N}| \leq \varepsilon} \right) \delta_N(dt), $$
where $\delta_N(dt) = \frac{N}{n(N)} \mathbf{1}_{\{ \frac{[Nt_0]}{N}-\frac{n(N)}{2N} \leq t < \frac{[Nt_0]}{N}+\frac{n(N)}{2N}\}} dt$.
Use again the formula of the density of a standard $\alpha$-stable random variable $f_{\alpha}(x) = \frac{1}{\pi} \int_{0}^{+\infty} e^{-|t|^{\alpha}} \cos(tx) dt $ and the inequality (\ref{probaepsilon}) to obtain the convergence $\Var \left((\mu_{\frac{[Nt]}{N}}-z_{[Nt],N})\mathbf{1}_{|X_{[Nt],N}| \leq \varepsilon} \right) \rightarrow \sigma_t^2$ and 
$$\sum_{k=[Nt_0] - \frac{n(N)}{2}}^{[Nt_0] + \frac{n(N)}{2}-1}\limits \Var(X_{k,N} \mathbf{1}_{|X_{k,N}| \leq \varepsilon})  \rightarrow \sigma_{t_0}^2.$$
We conclude the proof using Theorem 4.1 of \cite{Pet} \Box

\end{document}